\documentclass[11pt]{article}
\usepackage{amsmath,amssymb,epsf}
\usepackage[french]{babel}
\usepackage[applemac]{inputenc}
\advance \textwidth 3cm
\advance \textheight 3cm
\usepackage{amsfonts}
\usepackage{latexsym}
\newtheorem{theorem}{Th\'eor\`eme}
\newtheorem{lemme}{Lemme}
\newtheorem{corollaire}{Corollaire}
\newtheorem{definition}{D\'efinition}

\newtheorem{prop}{Proprit}
\newtheorem{remarque}{Remarque}
\newenvironment{preuve}[1]{\par\noindent\underline{Preuve #1} :\quad}%
{\unskip\nobreak\hfil\penalty50\hskip2em\null\nobreak\hfil%
$\Box$\parfillskip0pt\par\medskip}

\newcommand{\Tr}{\mathrm {Tr}}

\newcommand{\nn}{\mathbb N}
\newcommand{\rr}{\mathbb R}

\newcommand{\trace}{\mathrm{Tr}}

\newcommand{\demi}{\frac{1}{2}}

\title{Inversion des matrices de Toeplitz dont le symbole admet un z\'ero d'ordre fractionnaire positif, 
valeur propre minimale.}

\author{ Philippe Rambour\thanks{Universit\'{e} de Paris Sud,
      B\^atiment 425; F-91405
Orsay Cedex;
tel : 01 69 15 57 28 ; fax 01 69 15 60 19
      \mbox{e-mail : philippe.rambour@math.u-psud.fr}}
       \and Abdellatif Seghier\thanks{Universit\'{e} de Paris Sud,
        B\^atiment 425; F-91405
Orsay Cedex;
tel : 01 69 15 57 29 ; fax 01 69 15 72 34
       \mbox{ e-mail : abdelatif.seghier@math.u-psud.fr}}}
\date{}
\begin{document}
\maketitle
  \renewcommand{\abstractname}{Rsum}
     \begin{abstract}
\textbf{ Inversion des matrices de Toeplitz dont le symbole admet un z\'ero d'ordre fractionnaire positif, 
valeur propre minimale et op\'erateur diff\'erentiel assoc\'i\'e.}\\
     Cet article pr\'esente trois r\'esultats distincts.  Dans une premi\`ere partie nous donnons l'asymptotique 
     quand $N$ tend vers l'infini des coefficients des polynmes orthogonaux  de degr\'e $N$ associ\'es au poids 
     $\varphi_{\alpha}(\theta)=\vert 1- e^{i \theta} \vert ^{2\alpha} f_{1}(e^{i \theta})$, 
     o\`u $f_{1}$ est une fonction strictement positive suffisamment r\'eguli\`ere
     et  $\alpha> \frac{1}{2}, \quad 
     \alpha\in \mathbb R\setminus N$. Nous en d\'eduisons  l'asymptotique 
     des \'el\'ements de l'inverse de la matrice de Toeplitz $T_{N}(\varphi_{\alpha})$ au moyen d'un noyau int\'egral 
     $G_{\alpha}.$ 
     Nous prolongeons ensuite un r\'esultat de A. B\"ottcher et H. Windom relatif \`a l'asymptotique de la valeur 
     propre minimale des matrices de Toepliz de symbole $\varphi_{\alpha}$.  On sait que dans ce cas la plus petite valeur propre de cette matrice
   admet une asymptotique, quand $N$ tend vers l'nfini, de la forme $\frac{c_{\alpha}}{N^{2\alpha}} f_{1}(1)$. Pour $\alpha\in \mathbb N^*$
   A. B\"ottcher et H. Windom obtiennent une asymptotique de $c_{\alpha}$ quand $\alpha$ tend vers l'infini, et un encadrement de $c_{\alpha}$
   dans les autres cas. Nous obtenons ici le m\^eme type de r\'esultat, mais pour $\alpha$ r\'eel positif.      \end{abstract}
     \renewcommand{\abstractname}{Abstract}
          \begin{abstract}
\textbf{Inversion of Toeplitz matrices with singular symbol. Minimal eigenvalues. }\\
Three results are stated in this paper. The first one is devoted to the study of the orthogonal polynomial with respect of the weight $\varphi_{\alpha} (\theta)=\vert 1- e^{i \theta} \vert ^{2\alpha} f_{1}(e^{i \theta})$, with 
$\alpha> \demi$ and $\alpha \in \rr \setminus \nn $, and $f_{1}$ a regular function. We obtain an asymptotic expansion of the coefficients of these polynomials, and we deduce an asymptotic of the entries of 
$\left( T_{N} (\varphi_{\alpha})\right)^{-1}$ where $T_{N} (\varphi_{\alpha})$ is a Toeplitz matrix with symbol 
$\varphi_{\alpha}$. Then we extend a result of A. B\"ottcher and H. Widom result related to the minimal eigenvalue of the Toeplitz matrix $T_{N}(\varphi_{\alpha})$. For $N$ goes to the infinity it is  well known 
that this minimal eigenvalue admit as asymptotic $\frac{c_{\alpha}}{N^{2\alpha}} f_{1}(1)$. When $\alpha\in \nn$
the previous authors obtain an asymptotic of $c_{\alpha}$ for $\alpha$ going to the infinity, and they have the bounds of $c_{\alpha}$ for the other cases. Here we obtain the same type of results but for $\alpha$ a positive real.  
  \end{abstract}


\textbf{\large{Mots clef}}

\textbf{ Inversion des matrices de Toeplitz , valeur propre minimale.}
\section {Introduction}
Rappelons que si $f$ est une fonction de $L^1(\mathbb T)$
 on appelle matrice de Toeplitz d'ordre $N$ de symbole $f$ 
 et on note $T_N(f)$ la matrice $(N+1) \times (N+1)$ telle que 
 $\left(T_N (f)\right) _{k+1,l+1} =\hat f (l-k) \quad \forall \,k,l \quad 
 0\le k,l \le N$, en notant $\hat h (j)$ le coefficient de Fourier 
 d'ordre $j$ d'une fonction $h$ (voir \cite{Bo.3}). Ici nous nous intressons plus prcisment aux matrices de Toeplitz de symbole 
 $f=\vert 1 - \chi \vert ^{2 \alpha} f_1$ o $f_1$ est une fonction 
 de $L^1(\mathbb T)$ strictement positive sur le tore et o
 $\alpha$ est un rel strictement suprieur  $-\demi$. 
 
 Le but de cet article est, dans une premire partie, de complter les rsultats de \cite{RS04} et de \cite{RS10}. C'est  dire de donner une expression asymptotique des termes 
$\left(T_N (f)\right)^{-1}_{k,l}$ et des termes $\left(T_N (f)\right)^{-1}_{k,1}$ avec $k=[Nx]$ et $l=[Ny]$. Autrement dit nous nous intressons aux lments du "coeur" de l'inverse de la matrice
 $T_N (f)$ et au "coeur" de l'ensemble des coefficients du 
 polynme orthogonal de degr $N$ associ au symbole $f$.
Nous avons fourni ces rsultats pour $\alpha\in ]-\demi, \demi]$
dans \cite{RS10} et pour $\alpha $ entier positif dans \cite{RS04} 
(voir aussi, bien sr, \cite{W3}, et \cite{Bot}). Le r\'esultat pour $\alpha=1$ est d \`a Spitzer-Stone  (\cite{SpSt}) et aussi 
\`a Courant, Friedrichs, et Lewy (\cite{CFL}). Ici nous donnons 
l'expression de ces asymptotiques pour $\alpha$ non entier 
suprieur  $\demi$ ( voir les thormes \ref{predicteur}
 et \ref{inverse}). On remarque que ces asymptotiques sont donns par des noyaux dont l'expression est toujours la mme
dans le cas des coefficients $\left(T_N (f)\right)^{-1}_{k,1}$. 
De plus les noyaux $G_\alpha$ donnant l'asymptotique de $\left(T_N (f)\right)^{-1}_{k,l}$ sont de la mme forme pour tous les $\alpha$ dans $]0, + \infty[$.\\
 Un autre problme classique li aux matrices de Toeplitz 
est celui de leurs valeurs propres et particulirement de leurs 
valeurs propres extrmales ( \cite{GS},\cite{W2}, \cite{Wid}, \cite{W3}, \cite{GZ}) et notamment de 
la plus petite. Si $\lambda_\alpha $
dsigne la valeur propre minimale de $T_N (f)$ avec 
$f= \vert 1-\chi \vert ^{2\alpha} f_1$, avec $f_{1}$ fonction r\'eguli\`ere,  on sait que dans le cas o\`u l'exposant $\alpha$
est un entier naturel strictement positif
 (\cite{ BoW}, \cite{BoW2}, )  il existe une 
constante $c_\alpha $ telle que pour $N$ tendant vers l'infini 
$$\lambda_\alpha \sim f_1(1) \frac{c_\alpha}{N^{2\alpha}}.$$
Kac, Murder et Szeg\"{o} ont montr\'e dans \cite{KMS} que $c_{1}$ existe et vaut $=\pi^2$ , et Parter a obtenu 
dans \cite{Part} que l'existence de $c_{2}$ et $c_{2}=500,5467 $.
Dans \cite{Part2} et \cite{Part3} Parter a montrer l'existence de $c_{\alpha}$ dans le cas g\'en\'eral pour des exposants entiers. 
 C'est donc un problme ancien que d'\'evaluer ou au moins d'encadrer le rel $c_\alpha$
 et aussi de chercher un asymptotique de cette constante
 quand $\alpha $ tend vers l'infini. Dans \cite{BoW} B\"ottcher et Widom rpondent  cette question pour $\alpha$ entier.   Nous reprenons ici ce problme pour des exposants $\alpha$ qui sont des rels positifs
 non entiers (thorme \ref{valeurpropre}).
 Dans le th\'eor\`eme \ref{valeurpropre} nous montrons 
 l'existence des constantes $c_{\alpha}$ et nous  donnons une expression formelle de $c_\alpha$ pour tout $\alpha $ positif.
 Nous en d\'eduisons ensuite des encadrements des constantes 
 $c_\alpha$ pour les valeurs de $\alpha$ non entires 
 situes dans les intervalles $ ]0, \demi[$, $[\demi,1[$ et 
 $]1, + \infty[$ (corollaire \ref{encadrements}). Enfin nous retrouvons dans le cas o\`u $\alpha$ est un rel positif l'asymptotique de $c_\alpha$ 
 donn par B\"ottcher et Widom dans le cas entier (thorme 
 \ref{asymptotique}).
 \\ Dans le cas d'exposants r\'eels n\'egatifs, $\alpha \in ]-\demi,0[$ on pourra consulter \cite{BoG}

\section{Inverse des matrices de Toeplitz de symbole $ \vert 1- e^{i \theta}\vert ^{2 \alpha} f_{1}$.}
Nous consid\'erons maintenant les matrices de Toeplitz 
$T_{N}(\vert 1- e^{i \theta}\vert ^{2 \alpha} f_{1})$
o\`u $\alpha$ strictement sup\'erieur \`a $\frac{1}{2}$,
la fonction $f_{1}$ \'etant r\'eguli\`ere, appartenant \`a 
l'ensemble $A(\mathbb T, 3/2) $
avec $A(\mathbb T,\nu) =\{ c \in L^2(\mathbb T) / \sum_{n\in \mathbb Z} n^{\nu} \vert \hat c (n)\vert <\infty\}$,
si l'on appelle $\hat c(n)$ le coefficient de Fourier d'ordre $n$ d'une fonction $c$.
Pour de telles fonctions $f$ et $f_{1}$ nous savons (voir \cite {GZ}) qu'il existe des fonctions $g$ et $g_{1}$ dans 
$H_{2}^+ = \{c \in L^2(\mathbb T) / \hat c (n) = 0 \iff n<0\}$ telles que 
$\vert 1- e^{i \theta}\vert ^{2 \alpha} f_{1}=g \bar g$, 
$ f_{1} =g_{1} \bar {g_{1}}$ (avec $g = (1-e^{i \theta}) g_{1}$). 
Dans la suite nous notons  $\beta_{0}^{(\alpha)} = \widehat{(1/g)} (0)$ et $\chi(\theta)= e^{i \theta} $ .
 Le but de cette section est d'\'etablir compl\`etement les deux th\'eor\`emes qui suivent.
 Il faut noter que pour $\alpha$ entier naturel strictement sup\'erieur \`a 1 ces th\'eor\`emes ont \'et\'e compl\`etement \'etablis dans \cite{RS04},
 et dans le cas o\`u $\alpha$ est un r\'eel contenu dans l'intervalle $]\frac{-1}{2}, \frac{1}{2}]$ ils ont \'et\'e d\'emontr\'es dans \cite{RS10}. Ici
 nous en donnons les d\'emonstrations pour $\alpha>\frac{1}{2}$. On peut remarquer que dans le cas o\`u $\alpha$ est n\'egatif 
 le noyau intervenant dans l'expression des coefficients de l'inverse est diff\'erent de celui que nous donnons ici. Il est par contre identique 
 dans le cas o\`u $\alpha\in ]0, \frac{1}{2}]$ et pour $\alpha\in \mathbb N$.
 \begin{theorem} \label{predicteur}
        Si $ \alpha > \frac{1}{2} $, 
     nous obtenons pour $ 0 < x < 1,$
        $$ g_{1} (1) \Bigl(T_{N} (\vert 1
-\chi\vert^{2\alpha} f_1)\Bigr)_{[Nx]+1,1}^{-1}
      =   \overline{ \beta_{0}^{(\alpha)}} N^{\alpha-1}  \frac{1}{\Gamma (\alpha) } x^{\alpha-1}
                  (1-x)^\alpha +o(N^{\alpha-1}),$$
   uniform\'ement pour $x$ dans
 $[\delta_1,\delta_2] $,  $0<\delta_1<\delta_2 <1$.
  \end{theorem}
 
 \begin{theorem}\label{inverse}
                    Pour 
                    $0<x,y<1$ et pour  $\alpha > \frac{1}{2}$ nous obtenons 
             $$ f_1(1) \Bigl(T_{N}
         (\vert 1- \chi \vert ^{2 \alpha} f_{1})\Bigr)_{[Nx]+1,
         [Ny]+1}^{-1}
      = N^{2\alpha-1}\frac{1}{\Gamma^{2}(\alpha)} \Bigl (
      G_{\alpha}(x,y) \Bigr)
      +o(N^{2\alpha-1})$$
 uniform\'ement en $(x,y)$  pour $ 0 < \delta _{1} \leq x,y \leq \delta _{2} < 1$
  avec
$$ G_{\alpha} (x,y) = x^\alpha y^\alpha \int_{\max(x,y)} ^1
        \frac{(t-x)^{\alpha-1}(t-y)^{\alpha-1}}{t^{2\alpha}} dt.$$
             \end{theorem}
Rappelons que pour $\alpha=1$ cet \'enonc\'e est connu depuis tr\'es longtemps : c'est le th\'eor\`eme de Spitzer-Stone \cite{SpSt}. 

\begin{corollaire} \label{trace}
Si $\alpha >\frac{1}{2}$ 
on a, avec les m\^{e}mes hypoth\`eses que pour le th\'eor\`eme \ref{inverse}, 
$$f_{1}(1) \trace\left( \left( T_{N}( \vert 1- \chi \vert ^{2 \alpha} f_{1}\right)^{-1}\right) = N^{2\alpha} 
\frac{B(2\alpha,2\alpha)}{\Gamma ^2(\alpha) (2\alpha-1)} +o(N^{2\alpha}).$$
\end{corollaire}
Le r\'esultat pour $\alpha \in \mathbb N$ a \'et\'e obtenu par A. B\"{o}ttcher dans \cite{Bot}. On obtient 
$$ f_{1}(1) \trace\left( \left( T_{N}( \vert 1- \chi \vert ^{2 \alpha} f_{1}\right)^{-1}\right) = N^{2\alpha} 
\frac{(2\alpha-1)!(2\alpha-2)!}{(4\alpha-1)![(\alpha-1)!]^2}+o(N^{2\alpha}).
$$
 \section {Valeurs propres minimales.}
\begin{remarque}
Pour $\alpha > \frac{1}{2}$ nous pouvons prolonger la fonction $G_{\alpha}$ \`a $x=y$.
Ce n'est \'evidemment pas possible pour $\alpha\in [0, \frac{1}{2}]$.  Nous pouvons n\'eanmoins \'ecrire le th\'eor\`eme ( \ref {valeurpropre}),
en faisant un passage \`a la limite dans le cas o\`u $\alpha < \frac{1}{2}$. 
\end{remarque}
Rappelons tout d'abord la dfinition suivante

\begin{definition}
Si $f$ est une fonction de $L^2 ([0,1] \times[0,1] )$ et $n$ un 
entier naturel, nous noterons : 
$ \ast ^n f $ la fonction dfinie sur $[0,1] \times[0,1]$ par 
$$ 
\ast ^n f (x,y)=\int _{0}^1\int _{0}^1 f(x,x_1)\int_{0}^1f(x_{1},x_{2})
\cdots \int_{0}^1 f(x_{n-1},x_{n}) f(x_{n},y) dx_{n}dx_{n-1} 
\cdots  dx_{2}dx_{1} 
$$
\end{definition}

\begin{theorem} \label{valeurpropre}
Si $f_{1}\in A(\mathbb T, 3/2)$ et si 
$\lambda_{min} \left( T_{N}(\vert 1- \chi \vert ^{2\alpha} f_{1}) \right)$ d\'esigne la valeur propre minimale de $T_{N}(\vert 1- \chi \vert ^{2\alpha} f_{1})$
avec $ \lambda_{min} \left( T_{N}(\vert 1- \chi \vert ^{2\alpha} f_{1}) \right) \sim \frac{c_{\alpha}}{N^{2\alpha}} f_{1}(1)$ 
on a pour tout r\'eel $\alpha>0$  diffrent de $\demi$, 
$$\frac{f_{1}(1)} {N^{2\alpha}}\Bigl( \lim_{s\rightarrow + \infty } \left(\int_{0}^1 \ast ^s G_{\alpha} (t,t) dt \right) ^{1/s} \Bigr)= 
\frac{1}{c_{\alpha}} \left(1+o(1)\right).$$
\end{theorem}
\begin{corollaire}\label{encadrements}
En gardant les m\^emes notations et hypothses que pour le thorme prcdent on a les encadrements suivants pour la constante $c_{\alpha}$.
\begin{enumerate} 
\item
Si  $0<\alpha<\frac{1}{2}$ alors 
$$ \frac{  \Gamma(1-\alpha)}{\Gamma(1-2\alpha)}\frac{\Gamma(\alpha)}{K_{\alpha}}  \le c_{\alpha}\le 
\frac{ \Gamma(\alpha)  (4\alpha+1) \Gamma(1-a)}{\Gamma(1-2\alpha)}$$
avec 
$$K_{\alpha}= \left (\frac{1}{2\alpha} +\frac{\Gamma^2(2\alpha)}{\Gamma(4\alpha) } 
+\frac{\Gamma(1-2\alpha) \Gamma (\alpha)}{\Gamma(1-\alpha)}\right)$$
 \item
Si  $\frac{1}{2} < \alpha<1$ alors 
$$\frac{\Gamma(4\alpha) \Gamma^2(\alpha)(2\alpha-1)}{\Gamma^2(2\alpha) }\le c_{\alpha} 
\le \Gamma^2(\alpha)   \frac{ (2\alpha +1) (2\alpha+2) (2\alpha+3)}{2}.$$
\item
Si $ \alpha \in ]1, + \infty[ \setminus \mathbb N^*$ alors 
$$ \frac{\Gamma^2 (\alpha) \Gamma(4\alpha) } {\Gamma(2\alpha-1)\Gamma(2\alpha+1)}  \le c_{\alpha} \le
\frac{\Gamma^2(\alpha)}{2} (2\alpha-1) (4\alpha-1) (2^{4\alpha-1}) .$$
\end{enumerate}
\end{corollaire} 
Parmi les encadrements possibles de $c_{\alpha}$ nous avons essay\'e de donner les encadrements les meilleurs. 
Pour $\alpha$ entier sup\'erieur ou \'egal \`a $1$ B\"{o}ttcher et Widom obtiennent dans \cite{BoW} l`'encadrement 
$$\frac{\Gamma^2 (\alpha) \Gamma(4\alpha) } {\Gamma(2\alpha-1)\Gamma(2\alpha+1)}  \le c_{\alpha} \le
\frac{4\alpha+1}{2\alpha+1} \frac{\Gamma(4\alpha+1) \Gamma^2(\alpha+1)}{\Gamma^2(2\alpha+1)}.$$
Nous ne pouvons pas r\'ecup\'erer le majorant dans le cas o\`u $\alpha$ est un r\'eel positif non entier,
les arguments utilis\'es \'etant sp\'ecifiques aux entiers. N\'eanmoins notre majorant est d' ordre comparable. 

La valeur propre minimale est beaucoup plus difficile  valuer dans le cas 
$\alpha =\demi$. Nanmoins en utilisant le corollaire 3
 de \cite{RS10} que nous rappelons ici
 \begin{theorem} 
 Si $f_{1}$ est une fonction r\'eguli\`ere dans $\mathbb A(\mathbb T, \frac{3}{2})$ on a 
 $$ f_{1}(1) \Tr \left( T_{N} \left(\vert 1-\chi \vert f_{1}) \right) ^{-1} \right) = \frac{1}{\pi} N \ln N + o(N \ln N).$$
 \end{theorem}
nous obtenons imm\'ediatement la proprit
\begin{prop} \label{mi}
Si $\alpha = \demi$ on a 
$$ \lambda_{min}\left(T_{N}(\vert 1- \chi \vert ^{2\alpha} f_{1})\right) \ge \frac{\pi}{N \ln N}.$$
\end{prop}
Enfin de la mme manire que B\"{o}ttcher et Widom  dans \cite{RS10}nous obtenons l'asymptotique de $c_\alpha$ lorsque $\alpha$ tend vers l'infini.
\begin{theorem} \label{asymptotique}
Avec les hypoth\`eses du corollaire pr\'ec\'edent nous pouvons \'ecrire, pour $\alpha$ tendant vers plus l'infini 
$$  c_{\alpha}= \sqrt{8 \pi \alpha} \left( \frac{4 \alpha} {e} \right) ^{2\alpha}  \left( 1+O(1/\sqrt{\alpha})\right).$$
\end{theorem}
 \section{D\'emonstration du th\'eor\'eme \ref{predicteur}.
Cas o\`u $\alpha>\frac{1}{2}.$}
\subsection{Position du probl\`eme}
 Dans cette partie nous consid\'ererons un symbole de la forme 
 $f=\vert 1-\chi \vert ^{2\alpha} f_{1}$ avec $f_{1}$ une fonction r\'eguli\`ere qui s'\'ecrit
 $f_{1}= g_{1} \bar g_{1}$ et o $\alpha$ est un rel strictement 
 suprieur  $\demi$. dans la suite de la dmonstration nous supposerons,
 pour simplifier les notations, que $g_1(1)\in \mathbb R$. Pour de tels $\alpha$ nous posons 
 $g=(1-\chi)^\alpha g_1$ et $\frac{1}{g } = \displaystyle{\sum_{u=0} ^{+ \infty} \beta_{u}^{(\alpha)} \chi^u}$.
 Nous avons ainsi  $\beta_0^{(\alpha)}= \frac{1}{2 \pi} \int _{-\pi}^{\pi} \frac{1}{g(e^{i \theta})} d\theta $.
 Rappelons que les th\'eor\`emes \ref{predicteur} et \ref{inverse} ont \'et\'e \'etablis dans \cite{RS10}
 pour  $\alpha \in ]- \frac{1}{2},\frac{1}{2}] $ et dans \cite {RS04} pour $\alpha \in \mathbb N^*$. On pourra aussi consulter 
 \cite{RS12}.
 On peut noter que dans le cas particulier 
 $\alpha=1$ le th\'eor\`eme \ref{inverse} correspond au c\'el\`ebre th\'eor\`eme de Spitzer-Stone (voir \cite{SpSt}).
  Nous utiliserons pour cela d'une mani\`ere d\'eterminante la propri\'et\'e fondamentale des polyn\^{o}mes pr\'edicteurs
  dont les coefficients 
 sont obtenus en normalisant les termes de la premi\`ere colonne de l'inverse de$ T_{N}(f)$ par $\left(T_{N}(f) \right)_{1,1}^{1/2}$
(voir \cite{Ld}).\\
Rappelons ici cette proprit ainsi que la formule de Gohberg-Semencul \cite{GoSe}.
\begin{prop}\label{fondamentale}
Si $P_N$ dsigne le polyn\^ome prdicteur du symbole $h$ alors 
$$
 \forall s \quad -N\le s \le N \quad \widehat {h}(s) = \widehat{\left(\frac{1}{P_N}\right)} (s),
$$
si $\widehat {h}(s)$ dsigne le coefficient de Fourier d'indice $s$ de la fonction $h$.
On a alors  
\begin{equation}
T_N(h) = T_N \left(\frac{1}{P_N}  \right).
\end{equation}
\end{prop}
\begin{prop} (Goberg-Semencul)
Si $P_{N+1} =\displaystyle{\sum_{u=0}^{N+1} \beta_u \chi ^u }$
un polyn\^ome trigonomtrique de degr infrieur ou gal  $N+1$
on a , si $k\le l$ 
$$
T_N(\frac{1}{\vert P_{N+1}\vert ^2} )^{-1}_{k+1,l+1} 
=  \sum_{u=0}^k \bar \beta_{k-u} \beta_{l-u} -\sum_{N-l}^{N+k-l} \beta_u \bar \beta_{u+l-k}. 
$$
\end{prop}
 On remarque que la formule de Gohberg-Semencul et aussi la propri\'et\'e \ref{fondamentale}  
 permettent de calculer, en toute g\'en\'eralit\'e, les coefficients \\ $\left(T_{N}(f)\right)^{-1}_{h+1, l+1},$ $ 0 \le h\le N, \quad 0 \le l \le N$ 
 quand on conna\^{i}t les 
 coefficients $\left(T_{N}(f)\right)^{-1}_{k+1, 1} $ $0 \le k \le N.$
 La d\'emarche naturelle est donc de d'obtenir d'abord le th\'eor\`eme \ref{predicteur} et d'en d\'eduire le th\'eor\`eme \ref{inverse}.\\
Dans le cas o\`u $\alpha>\demi$ la propri\'et\'e fondamentale des polyn\^{o}mes pr\'edicteurs et la formule de 
 r\'ecursion \ref{recurcul} ci-dessous permettent de d\'eterminer les termes 
 $\left(T_{N}(\vert 1-\chi \vert ^{2\alpha} f_{1})\right)^{-1}_{k+1, 1}$ pour \mbox{$0 \le k \le N$}
en fonction des termes $\left(T_{N}(\vert 1-\chi \vert ^{2\alpha-2} f_{1})\right)^{-1}_{u+1, 1} \quad 0 \le u \le N.$
 Il est donc naturel d'utiliser un raisonnement par r\'ecurrence pour obtenir les coefficients 
 $\left(T_{N}(\vert 1-\chi \vert ^{2\alpha} f_{1})\right)^{-1}_{u+1, 1}$, $ 0 \le u \le N,$ pour tous les r\'eels $\alpha>  \frac{1}{2}$,
 $\alpha \not \in \mathbb N$ 
 \`a partir des \'el\'ements
  $\left(T_{N}(\vert 1-\chi \vert ^{2\alpha-2} f_{1})\right)^{-1}_{k+1, 1},$ $ 0 \le k \le N$. Dans le cas o\`u
   $- \frac{1}{2}<\alpha \le\frac{1}{2}$ ces coefficients sont
connus explicitement ( voir \cite {RS10}), nous pouvons donc initialiser sans probl\`eme la r\'ecurrence. Plus pr\'ecis\'ement 
  pour un entier positif $p$ et si $\alpha$ est un r\'eel appartenant \`a l'ensemble
  $]p+ \frac{1}{2}, p + \frac{3}{2}]\setminus \{p\}$ 
  nous allons d\'eduire, \`a l'aide de la formule \ref{recurcul}, les coefficients de la premi\`ere colonne de $\left(T_{N}(\vert 1-\chi \vert ^{2\alpha} f_{1})\right)^{-1}$
  de ceux de   $\left(T_{N}(\vert 1-\chi \vert ^{2\alpha-2} f_{1})\right)^{-1}$ que nous supposerons connus.
  \subsection{Formules pr\'eliminaires.}
  Citons tout d'abord le rsultat suivant tabli, comme nous l'avons dj 
  dit, dans \cite{RS10}.
 \begin{theorem} \label{Theo1}
	 Si $\frac{-1}{2}<\alpha< \frac{1}{2}$
      et $ 0<x<1$   
      $$  g_1(1) \left(T_{N} (\vert 1-\chi \vert ^{2\alpha} f_{1})
      \right)^{-1}_{[Nx]+1,1}= \overline{\beta_{0}^{(\alpha)} }\frac{ N^{\alpha-1}}{\Gamma (\alpha)} x^{\alpha-1} (1-x)^\alpha
      + 
      o(N^{\alpha-1})$$
      uniform\'{e}ment par rapport \`{a} $x$, sur tout intervalle 
      $[\delta_1,\delta_2],$ $0<\delta_1<\delta_2<1.$      
       \end{theorem} 
       \begin {theorem} \label{Theo1bis}
Avec les m\^emes hypoth\`eses que dans le th\'eor\`eme pr\'ec\'edent nous avons, pour tout
 $kÊ\in \nn$ tel que 
$\displaystyle{\lim_{N\rightarrow \infty}\frac {k}{ N}=0}$,
les relations suivantes
\begin{enumerate}
\item
\begin{equation}
\left(T_N(\vert 1-\chi \vert ^{2\alpha} f_{1})\right)^{-1}_{k+1,1}= \overline{ \beta_{0}^{(\alpha)}}\left(\beta^{(\alpha)}_k-
\frac{\alpha^2} {N} \beta_k^{(\alpha+1)}\right)\left(1+o(1)\right),~~ N\rightarrow\infty
\end{equation}
\item
\begin{equation}
\left(T_{N}(\vert 1-\chi \vert ^{2\alpha} f_{1})\right)^{-1}_{N+1-k,1}=\overline{ \beta_{0}^{(\alpha)}}\left(
\frac {g_ (1))}
{\bar g_1(1)} \beta_k^{(\alpha+1)}\right) \frac {\alpha}
{N} \left(1+o(1)\right) ,~N\rightarrow \infty,
\end{equation}
\end {enumerate}
 les deux formules \'etant uniformes en $k$ pour $k \in [0, [N \epsilon]]$, $\epsilon$ \'etant un entier suffisamment petit.
\end{theorem}
  (On pourra consulter \cite{KaRS} pour la d\'emonstration de ce dernier th\'eor\`eme).

                 Dans la suite de la d\'{e}monstration nous allons 
     utiliser la formule de r\'ecursion suivante permettant de calculer pour tout polyn\^{o}me de degr\'e $N+1$ les termes 
     $\left(T_{N, \vert 1-\chi\vert ^2/ 
     \vert P_{N+1}\vert ^2}\right)^{-1}_{k+1,1}$ en fonction des coefficients de $P_{N+1}$. Cette formule a \'{e}t\'{e}
     \'{e}tablie dans \cite{RS04}. Comme nous l'avons montr\'{e} 
     dans ce m\^eme article, elle permet de calculer le polyn\^ome 
     pr\'{e}dicteur de $\vert 1-\chi\vert ^2 h$ en fonction de celui 
     de $h$.  Si $P_{N+1} = \sum_{u=0}^{N+1} \beta_u \chi^u$ cette formule
     s'\'{e}crit : 
     \begin{equation}  \label{recurcul}
      \left ( T_{N, \vert 1-\chi\vert ^2/ \vert
     P_{N+1}\vert ^2}\right)^{-1}_{k+1,1} = \bar \beta_0 \sum_{u=0}^k \bar
     \beta_{k-u} + \frac{1}{N+1+A(P_{N+1})}
   \bar \beta_0 A_{N,k}
    \end{equation}
    avec :
     \begin{enumerate}
     \item
    $ A_{N,k}= \left(\tilde{Q}'_{2,k}(1)   
     {P}_{N+1}(1)/ \bar P_{N+1}(1)  - (Q'_{1,k}) (1) + \displaystyle 
     {\sum_{u=0}^k} 
     \beta_u\right),$
     \item 
          $A(P_{N+1}) = -2 \Re (\bar{P'}_{N+1}(1) P_{N+1}(1))/
     \vert P_{N+1}(1) \vert^2$
     \item
    $ \tilde{Q}_{2,k}(r) =\displaystyle{ \sum_{u=N+2-k}^{N+1} \bar\beta_u
    r^{u-(N+1)+k}}  \quad 
    Q_{1,k}(r) =\displaystyle
    {\sum_{u=0}^k \beta_u r^{k-u+2}}.$
    \end{enumerate} 
   
Enfin nous noterons par 
$ \frac {1}{N+2+ A (P_{N+1,\alpha})}$ o\`{u} $P_{N+1,\alpha}$
 le polyn\^ome pr\'{e}dicteur de
$\vert 1-\chi\vert ^{2\alpha} f_1$, 
 et nous poserons :
$$ P_{N+1,\alpha} = \sum_ {u=0}^{N+1} \gamma_{u}^{(\alpha)} \chi^u.$$ 
Il est facile de se convaincre que pour pouvoir utiliser la formule (\ref{recurcul}) pour passer de la connaissance de 
$\left (T_{N} (\vert 1-\chi\vert ^{2\alpha} f_{1}) \right)^{-1}_{k+1,1}$ \`a celle de 
$\left (T_{N} (\vert 1-\chi\vert ^{2\alpha+2} f_{1}) \right)^{-1}_{k+1,1}$ nous
avons besoin de conna\^{i}tre les quantit\'es 
\begin{itemize}
\item
$$ \left (T_{N} (\vert 1-\chi\vert ^{2\alpha} f_{1}) \right)^{-1}_{k+1,1} \quad \mathrm{quand} \quad
 \displaystyle{\lim _{N\rightarrow + \infty}  k/N =0}$$
 \item
$$ \left (T_{N} (\vert 1-\chi\vert ^{2\alpha} f_{1}) \right)^{-1}_{k+1,1} \quad \mathrm{quand}\quad
 \displaystyle{\lim _{N\rightarrow + \infty}  k/N =x,\, 0<x<1}$$
 \item
 $$ \left (T_{N} (\vert 1-\chi\vert ^{2\alpha} f_{1}) \right)^{-1}_{N-k+1,1} \quad \mathrm{quand}\quad
 \displaystyle{\lim_{N\rightarrow + \infty}  k/N =1}$$
 \item
$$ P_{N+1,\alpha} (1) \quad \mathrm{et} \quad P'_{N+1,\alpha} (1).$$
 \end{itemize}
 Pour $\alpha\in ]- \frac{1}{2}, \frac{1}{2}]$ nous connaissons les trois premi\`eres quantit\'es. Il nous faut obtenir les deux derni\`eres.  Ces valeurs peuvent facilement se d\'eduire des r\'esultats par un passage \`a la limite dans le th\'eor\`eme
 1.4 de l'article \cite{ML3}. 
 
 Nous sommes finalement conduits  \`a \'etablir l'hypoth\`ese de r\'ecurrence suivante, qui doit \^{e}tre v\'erifi\'ee pour tout 
 entier naturel $p$. \\
 \textbf{Hypoth\`ese de r\'ecurrence}
 Pour tout entier $p \in \mathbb N$ et tout r\'eel $\alpha \in ]p- \frac{1}{2}, p+ \frac{1}{2}]\setminus\{p\}$ nous avons les asymptotiques suivantes 
  \begin{enumerate}
 \item 
$$ \left (T_{N} (\vert 1-\chi\vert ^{2\alpha} f_{1}) \right)^{-1}_{k+1,1}  =  \overline {\beta_{0}^{(\alpha)}} \beta^{(\alpha)}_{k} \left(1+o(1)\right)\quad \mathrm{pour} \quad
 \displaystyle{\lim _{N\rightarrow + \infty}  k/N =0}$$
\item
\centerline{Si
$ \displaystyle{\lim _{N\rightarrow + \infty}  k/N =x,\, 0<x<1}$
 nous avons} 
$$ \left (T_{N} (\vert 1-\chi\vert ^{2\alpha} f_{1}) \right)^{-1}_{k+1,1} =\overline{ \beta_{0}^{(\alpha)}} \frac{N^{\alpha-1}} {\Gamma (\alpha) g_1(1)}
 x^{\alpha-1} (1-x)^\alpha +o(N^{\alpha-1})$$
  
\item
$$ \left (T_{N} (\vert 1-\chi\vert ^{2\alpha} f_{1}) \right)^{-1}_{N-k+1,1}= O \left(\max (N^{\alpha-1}), \frac{1}{N}\right)
\quad \mathrm{pour}\quad
 \displaystyle{\lim_{N\rightarrow + \infty}  k/N =0}$$
 \item 
 $$ P_{N+1,\alpha} (1) = \frac{N^\alpha}{g_{1}(1)} \frac{\Gamma (\alpha+1)}{\Gamma(2\alpha+1)}+o(N^\alpha),
 \quad 
 P'_{N+1,\alpha} (1) = \frac{N^{\alpha+1}}{g_{1}(1)} \frac{\Gamma^2 (\alpha+1)}{\Gamma(2\alpha+2)\Gamma (\alpha)}+o(N^{\alpha+1}).
 $$

 \end{enumerate}
 Pour $p=0$ ( c'est \`a dire $\alpha \in ]-\demi, \demi[$) ces points sont \'etablis : ce sont les th\'eor\`emes \ref{Theo1} et 
 \ref{Theo1bis} pour les points 1. 2. et 3. et un passage \`a la limite dans 
 \cite{ML3}. Il nous faut maintenant les \'etablir pour l'entier $p+1$ en les supposant vrais pour $p$. Cependant  
 certains calculs demandent un traitement diff\'erent selon que le r\'eel $\alpha$ est positif ou n\'egatif. Nous serons donc amen\'es \`a distinguer parfois ces deux cas dans nos d\'emonstrations.  Remarquons enfin que l'uniformit\'e de l'expression asymptotique passe de
$\alpha$ \`a $\alpha' =\alpha+1$  gr\^ace \`a la relation  (\ref{recurcul}).
  \subsection{ Calcul de l'asymptotique de $\left (T_{N} (\vert 1-\chi\vert ^{2\alpha} f_{1}) \right)^{-1}_{k+1,1} \,\mathrm{quand} \,
 \displaystyle{\lim _{N\rightarrow + \infty}  k/N =0}$}
Dans un premier temps il est facile de v\'erifier que  l'hypoth\`ese de r\'ecurrence, jointe \`a la formule \ref{recurcul} et 
\`a la normalisation du polyn\^ome pr\'edicteur, donne, pour tout  $\alpha >0$, 
$\gamma_0^{(\alpha+1)}=\gamma_0^{(\alpha)} =\beta_{0}^{(\alpha)}+ O(\frac{1}{N})$.  \\
Si maintenant $\frac{k}{N} \rightarrow 0$ quand $N$ tend vers l'infini, la formule  (\ref{recurcul}) 
devient
 $$\left( T_N ( \vert 1-\chi \vert ^2 / \vert P_{N+1} \vert) \right)^{-1} _{k+1,1} =
 \sum_{u=0}^k \gamma_u ^{(\alpha)} +\frac {2\alpha +1}{N} \left( -(k+1) \sum_{u=0}^k
  \gamma_u ^{(\alpha)} + \sum_{u=0}^k u\gamma_u ^{(\alpha)}\right) \left (1+o(1)\right).$$
  si $k$ prend une valeur fixe suffisamment petite ($k=0,1,2\cdots$), le r\'esultat est obtenu imm\'ediatement pour $\alpha+1$.
  Si ce n'est pas le cas nous pouvons remarquer que \\
  $ \displaystyle{ \sum_{u=0} ^k  \beta _u ^{(\alpha)} =\beta_k ^{(\alpha+1)}}$ et 
   $ \displaystyle{\sum_{u=0} ^k u \beta _u ^{(\alpha)} 
\sim\alpha \beta _{k-1} ^{(\alpha+2)}\sim \alpha \beta_k ^{(\alpha+1)}  \frac{k} {(\alpha+1)}}$.
 Nous pouvons alors \'ecrire 
  \begin{eqnarray*}
 & &  \left( T_N ( \vert 1-\chi \vert ^2  / \vert P_{N+1} \vert) \right)^{-1} _{k+1,1} =
  \\ &=& \overline{\beta_0^{(\alpha)}} \Bigl(
  \beta_k ^{(\alpha+1)} +\frac {2\alpha +1}{N} \left( -\ (k+1)  \beta_k ^{(\alpha+1)}
  +\alpha  \beta_k ^{(\alpha+1)}  \frac{k} {\alpha+1}\right) \Bigr)\left( 1+o(1)\right).
  \end{eqnarray*}
  C'est \`a dire que quand $\frac{k}{N} \rightarrow 0$ quand $N$ tend vers l'infini nous avons 
  $$\left( T_N ( \vert 1-\chi \vert ^2 / \vert P_{N+1} \vert) \right)^{-1} _{k+1,1} = 
 \overline{\beta_0^{(\alpha)}} \beta_k ^{(\alpha+1)}  \left( 1+o(1)\right),$$
ce qui s`'\'ecrit encore 
  $$\left( T_N ( \vert 1-\chi \vert ^{2(\alpha+1)} f_1 \right)^{-1} _{k+1,1} = 
\overline{\beta_0^{(\alpha+1)}}  \beta_k ^{(\alpha+1)}  \left( 1+o(1)\right).$$
Nous avons donc d\'emontr\'e le point 1. de la r\'ecurrence.
    \subsection{ Calcul de l'asymptotique de $\left (T_{N} (\vert 1-\chi\vert ^{2\alpha} f_{1}) \right)^{-1}_{N-k+1,1} \quad \mathrm{quand} \quad \\
 \displaystyle{\lim _{N\rightarrow + \infty}  k/N =1}$}
 L'hypoth\`ese de r\'ecurrence et la formule de r\'ecursion (\ref{recurcul}) permettent d'\'ecrire, en utilisant l'hypothse $g_1(1) \in \mathbb R$ 
\begin{eqnarray*} 
& & \left( T_N ( \vert 1-\chi \vert ^2 / \vert P_{N+1} \vert) \right)^{-1} _{k+1,1} 
  \sim  P_{N+1}(1) +\\
  &+&  \frac{2\alpha+1}{N} \left( P'_{N+1}(1) 
  -o(N) P_{N+1} (1) + \frac{ g_1(1)}{\bar g_1(1)}\left( \overline{P'_{N+1} (1)}
    -N \overline{P_{N+1}(1)} \right)\right) \\
  &\sim& O(N^\alpha) = O \left( (\max ( N^{\alpha'}, \frac{1}{N})\right). 
\end{eqnarray*}
L'\'egalit\'e $O(N^\alpha) = O \left( (\max ( N^{\alpha'-1}, \frac{1}{N})\right)$ \'etant une cons\'equence de  $p\ge0$ nous obtenons ainsi 
l'\'egalit\'e annonc\'ee, c'est \`a dire le point 3. de l'hypoth\`ese de r\'ecurrence.  Dans les paragraphes qui suivent nous allons maintenant
\'etablir le point 2. de cette hypoth\`ese.\\
 \subsection{ Calcul de l'asymptotique de $\left (T_{N} (\vert 1-\chi\vert ^{2\alpha} f_{1}) \right)^{-1}_{N-k+1,1} \quad \mathrm{quand} \\\
 \displaystyle{\lim _{N\rightarrow + \infty}  k/N = x}, 0<x<1.$}
 Dans ce cas nous pouvons \'ecrire la formule \ref{recurcul} de la mani\`ere suivante 
 \begin{align} \label{recurla}
&  \left( T_N ( \vert 1-\chi \vert ^2 / \vert P_{N+1} \vert) \right)^{-1} _{k+1,1} =\\
&= \bar \gamma^{(\alpha)}_0 \left( \sum_{u=0} ^k \gamma^{(\alpha)}_{k-u} + \frac{1}{N+1+A(P_{N+1})} \times \right.\\
& \times \left.\left( \sum _{u=N+2-k}^{N1}\overline{ \gamma^{(\alpha)}_u} \left(u-(N+1)+k\right) 
 \frac{P_{N+1}(1)}{\overline{P_{N+1}(1)}} - \sum_{u=0}^k \gamma^{(\alpha)}_u (k-u+1)\right)\right) \left(1+O(1/N)\right)\nonumber
\end{align}
\textbf {Calcul de la somme  $\displaystyle{\sum_{u=0} ^k \gamma^{(\alpha)}_{k-u}}$.}\\
 \textbf{Cas $\alpha<0$}.\\
 Ecrivons $ \displaystyle{\sum_{u=0} ^k \gamma^{(\alpha)}_{k-u}}= S_{1}+S_{2}$ avec $S_{1}=\displaystyle{ \sum_{u=0}^{[N \epsilon] }\gamma_{u}^{(\alpha)}}$ et 
 $S_{2}= \displaystyle{\sum_{[N \epsilon]+1}^k \gamma_{u}^{(\alpha)}}$, avec $\epsilon$ un r\'eel positif qui tend vers z\'ero.
 On trouve 
 $$ S_{1}=\frac{1}{\Gamma (\alpha+1) g_{1}(1) } [N \epsilon]^\alpha+ o(N^\alpha),$$
 et 
 \begin{eqnarray*}
  S_{2}&=& \frac{1}{\Gamma (\alpha) g_{1}(1) } N^\alpha \left ( \int_{[N \epsilon]+1)/N}^x t^{\alpha-1} 
 \left( (1-t)^\alpha-1\right) dt + \int_{[N \epsilon]+1)/N}^x t^{\alpha-1}\right)+ o(N^\alpha)\\
 &=& \frac{1}{\Gamma (\alpha) g_{1}(1) } N^\alpha \left ( \int_{[N \epsilon]+1)/N}^x t^{\alpha-1} 
 \left( (1-t)^\alpha-1\right) dt + \frac{x^\alpha}{\alpha} - \frac{[N \epsilon]^\alpha}{\alpha}\right)+ o(N^\alpha).
 \end{eqnarray*}
 Ce qui donne finalement 
 $$ \sum_{u=0} ^k \gamma^{(\alpha)}_{k-u} =\frac{1}{\Gamma (\alpha) g_{1}(1) } N^\alpha \left ( \int_{[N \epsilon]+1)/N}^x t^{\alpha-1} 
 \left( (1-t)^\alpha-1\right) dt + \frac{x^\alpha}{\alpha}\right) + o(N^\alpha).$$
 \textbf{Cas $\alpha>0$}.\\
La formule d'Euler Mac-Laurin donne imm\'ediatement 
 $$ \sum_{u=0} ^k \gamma^{(\alpha)}_{k-u} =
  \frac{1}{\Gamma (\alpha) g_{1}(1) } N^\alpha  \int_{0}^x t^{\alpha-1} (1-t)^\alpha dt + o(N^\alpha).$$
 \textbf {Calcul de la somme $\displaystyle{ \sum_{u=N+2-k}^{N+1} u \gamma_{u}^{(\alpha)} }.$}\\
L\`a aussi la formule d'Euler Mac-Laurin donne imm\'ediatement dans tous les cas 
$$  \sum_{u=N+2-k}^{N+1} u \gamma_{u}^{(\alpha)} = \frac{N^{\alpha+1}}{\Gamma(\alpha) g_{1}(1)} 
\int_{1-x} ^1 t^\alpha (1-t)^\alpha dt + o(N^{\alpha+1}).$$
\textbf {Calcul de la somme $\displaystyle{ \sum_{u=N+2-k}^{N+1}  \gamma_{u}^{(\alpha)} }.$}\\
 On obtient de m\^{e}me, en utilisant les points 1 et 3 de l'hypoth\`ese de r\'ecurrence 
 $$ \sum_{u=N+2-k}^{N+1}  \gamma_{u}^{(\alpha)} (k-N-1) =
 N^\alpha (k-N-1) \int_{1-x} ^x t^{\alpha-1} (1-t)^\alpha dt +o(N^{\alpha+1}) +o(N^{\alpha+1}).$$
 Le dernier point de l'hypoth\`ese de r\'ecurrence permet d'autre part d'affirmer que, nous avons, dans tous les cas 
 $$ A(P_{N+1,\alpha}) = - \frac{2\alpha}{2\alpha+1} N +o(N).$$
 \textbf {Fin de la d\'emonstration point 2.}
 Nous avons obtenu :
  \begin{itemize} 
  \item
  si $\alpha<0$  
 $$\left( T_{N}\left( \frac{\vert 1-\chi\vert ^{2(\alpha+1)}} {\vert P_{N+1}\vert^2}\right)\right)^{-1}_{k+1,1}
 =N^{\alpha} \frac{\bar \beta_{0}^{\alpha}}{\Gamma(\alpha) g_{1}(1) } k_{\alpha,1}(x) +o(N^{\alpha-1}),$$
 avec 
\begin{eqnarray*}
  k _{\alpha,1}(x) &= &
 \int _0 ^x t^{\alpha-1}\left( (1-t)^\alpha-1\right) dt  + \frac{x^\alpha}{\alpha}+\\ 
	    &+& (2\alpha+1)  \left( \int_{1-x} ^1 t^{\alpha-1} (1-t)^\alpha (t+x-1) dt  + \int _0^x t^\alpha (1-t)^\alpha dt \right.
\\
& -& \left. \frac{x^{\alpha+1}}{\alpha}
	 - x \int_0^x
	 t^{\alpha-1} \left((1-t)^\alpha-1 \right) dt \right),
\end{eqnarray*}

 \item 
 si $\alpha>0$  
 $$\left( T_{N}\left( \frac{\vert 1-\chi\vert ^{2(\alpha+1)}} {\vert P_{N+1}\vert^2}\right)\right)^{-1}_{k+1,1}
 =N^{\alpha} \frac{\bar \beta_{0}^{\alpha}}{\Gamma(\alpha) g_{1}(1) } k_{\alpha,2}(x) +o(N^{\alpha-1}),$$
 avec 
 \begin{eqnarray*}
  k _{\alpha,2}(x) &= &
 \int _0 ^x t^{\alpha-1} (1-t)^\alpha dt  +
	    + (2\alpha+1)  \left( \int_{1-x} ^1 t^{\alpha-1} (1-t)^\alpha (t+x-1) dt \right.
\\
& & \left. + \int_0^x t^{\alpha} (1-t)^\alpha dt 
	 - x \int_0^x
	 t^{\alpha-1} (1-t)^\alpha  dt \right).
\end{eqnarray*}
\end{itemize}
On v\'erifie par un calcul imm\'ediat que 
$$  \frac{d^2}{dx^2} k_{\alpha,1}(x)=\frac{d^2}{dx^2} k_{\alpha,2}(x) 
=\frac{1}{\alpha}\frac{d^2}{dx^2} \left( x^\alpha(1-x)^{\alpha+1}\right). $$
Puisque le d\'eveloppement en s\'erie enti\`ere de ces trois fonctions est une somme de termes en $t^{\alpha+j}, j \in \mathbb N^*$
ces fonctions ne peuvent diff\'erer d'un polyn\^ome de degr\'e un. Elles sont donc \'egales, ce qui ach\`eve de prouver le point 2 de la r\'ecurrence.
\subsection{ Calcul de $ P_{N+1,\alpha}(1) $ et de $P'_{N+1,\alpha}(1) $}
Ces quantit\'es sont \'etablis pour $\alpha \in ]-Ê\demi, \demi[$. Le calcul dans le cas g\'en\'eral est une simple application des points 1. 2. 3. 
de la r\'ecurrence et de la formule d'Euler et Mac-Laurin.

\section {D\'emonstration du th\'eor\`eme \ref{inverse} et du corollaire 
\ref{trace}}
\subsection{D\'emonstration du th\'eor\`eme \ref{inverse}}
Le thorme a t dmontr pour $\alpha \in ]-\demi, \demi]$
dans \cite{RS10}, et pour le cas o $\alpha$ est un entier dans 
$\mathbb N$ dans \cite{RS04}. Ici nous le dmontrons pour 
$\alpha \in ]p-\demi, p+\demi[$.Pour cela nous allons bien s\^ur 
utiliser les rsultats du thorme \ref{predicteur}.\\
 La formule de Gohberg-Semencul donne, pour $k\le l$, 
$$\left( T_{N}(\vert 1-\chi \vert ^{2\alpha-2} f_{1})\right)_{k+1,l+1}^{-1} = \sum _{u=0}^k \bar \gamma_{k-u}^{(\alpha)} \gamma_{l-u}^{(\alpha)} 
- \sum_{u=N-l}^{N+k-l}  \gamma_{v}^{(\alpha)}
\bar \gamma_{v+l-k} ^{(\alpha)},$$
avec toujours, si $k=[Nx], 0<x<1$ 
$$\gamma_{u}^{(\alpha)} = \frac{N^{\alpha-1}}{ \Gamma(\alpha) g_{1}(1) } x^{\alpha-1} (1-x)^\alpha
 + o(N^{\alpha-1}).$$
Posons $k=[Nx], l=[Ny], 0<x\le y<1$.
On peut dans un premier temps \'ecrire, en posant $k_{0}=[N \epsilon],\, \epsilon>0$ 
$$
\sum _{u=0}^k \bar \gamma_{k-u}^{(\alpha)}
 \gamma_{l-u}^{(\alpha )}=\sum _{u=0}^{k-k_{0}} \bar \gamma_{k-u}^{(\alpha)} \gamma_{l-u}^{(\alpha)}+
\sum _{u=k-k_{0}+1}^k \bar \gamma_{k-u}^{(\alpha)} \gamma_{l-u}^{(\alpha)}.
$$
On a clairement 
$$ \sum _{u=k-k_{0}+1}^k \bar \gamma_{k-u}^
{(\alpha)} \gamma_{l-u}^{(\alpha)}
= o(N^{2\alpha-1}).$$
D'autre part
\begin{eqnarray*}
&& \sum _ {u=0}^{k-k_{0}} \bar 
\gamma_{k-u}^{(\alpha)} \gamma_{l-u}^{(\alpha)}
 =\\
 &=& \frac{N^{2\alpha-2}} {\Gamma ^2 (\alpha) f_{1}(1)}  
 \sum _ {u=0}^{k-k_{0}} 
 \Bigl( \left(\frac{k-u} {N} \right)^{\alpha-1} \left ( 1- \frac{k-u}{N} \right)^{\alpha} 
\left(\frac{l-u} {N} \right)^{\alpha-1} \left ( 1- \frac{l-u}{N} \right)^{\alpha} +o(1) \Bigr )\\
&=& \frac{N^{2\alpha-1}} {\Gamma ^2 (\alpha) f_{1}(1)}  
\int_{0}^x (x-t)^{\alpha-1} (1-x+t)^{\alpha} (y-t)^{\alpha-1} (1-y+t)^\alpha dt + o(N^{2\alpha-1}).
\end{eqnarray*}
Si $J_{1}= \int_{0}^x (x-t)^{\alpha-1} (1-x+t)^{\alpha} (y-t)^{\alpha-1} (1-y+t)^\alpha dt$, nous pouvons \'ecrire,
en posant $u=\frac{t}{x},$ puis $w=1-u$, on a
$$ J_1 = x^\alpha \int_{0}^1 (1-u)^{\alpha-1} \left(1-x(1-u)\right) (y-ux)^{\alpha-1} (1-y+xt)^\alpha dt $$
$$J_1= x^\alpha \int_{0}^1 w^{\alpha-1} (1-xw)^\alpha \left(y-(1-w) x \right)^{\alpha-1}
\left(1-y+x(1-w) \right)^\alpha dw.$$
Enfin en posant $v = wx$ il vient
$$ J_1 = \int_{0}^x v^{\alpha-1} (1-v)^\alpha (y-x+v)^{\alpha-1}(1-y+x-v)^\alpha dv.$$
Nous avons de m\^{e}me  
$$
 \sum_{u=N-l}^{N+k-l}  \gamma_{v}^{(\alpha )} \bar \gamma_{v+l-k} ^{(\alpha)}
 =  \sum_{v=0} ^{k}   \gamma_{v+N-l}^{(\alpha)} 
\bar\gamma_{v+N-k} ^{(\alpha)} 
= \frac{N^{2\alpha-2 }}{Ê\Gamma ^2(\alpha) f_1(1)} S_{\alpha,k,N}
 $$ avec 
 \begin{eqnarray*}
S_{\alpha,k,N}&=& \Bigl(\sum_{v=0} ^{k} \left( \frac{v+N-k} {N} \right)^{\alpha -1}
\left(1- \frac{v+N-k} {N} \right)^{\alpha }
\left( \frac{v+N-l} {N} \right)^{\alpha -1}
\left(1- \frac{v+N-l} {N} \right)^{\alpha}+o(1)\Bigr) \\
&=& \frac{N^{2\alpha-1} }{Ê\Gamma ^2(\alpha) f_1(1)}
\int_0^x (1+t-x)^{\alpha-1} (x-t)^\alpha (1+t-y)^{\alpha-1}
(y-t)^{\alpha} dt+o(N^{2\alpha-1})
\end{eqnarray*}
Si l'on pose 
$$ J' =\int_0^x (1+t-x)^{\alpha-1} (x-t)^\alpha (1+t-y)^{\alpha-1}
(y-t)^{\alpha} dt$$
le changement de variables $u =\frac{t}{x}$ 
donne 
$$J' =x^{\alpha+1} \int_0^1 \left( 1- x(1-u)\right)^{\alpha-1} (1-u)^\alpha
(1-y+ux)^{\alpha-1} (y-ux)^{\alpha} du.$$
 L'on pose alors successivement 
 $ w= 1-u $ et $v= xw$. Cela donne d'abord 
 $$ J' =x^{\alpha+1} \int_0^1 \left( 1- xw\right)^{\alpha-1} w^\alpha
(1-y+x-wx)^{\alpha-1} (y-x+wx)^{\alpha} du,$$
puis
$$J' =\int_0^x \left( 1- v \right)^{\alpha-1} v^\alpha
(1-y+x-v)^{\alpha-1} (y-x+v)^{\alpha} dv.$$
Nous devons maintenant valuer $J-J'$. Nous obtenons tous calculs faits
\begin{eqnarray*}
 J_1-J' &=& \int _0^x v^{\alpha -1} (1-v)^{\alpha-1 } (y-x+v)^{\alpha -1} 
 (1-y+x-v)^{\alpha-1} (1-y+x-2v) dv\\
&=&\int_0 ^x (v-vy+vx-v^2)^{\alpha-1} (y-x+v-vy+vx-v^2)^{\alpha -1}
(1-y+x-2v) dv
\end{eqnarray*}
Effectuons maintenant le changement de variables 
$v_1 = v-vy +vx-v^2$. L'intgrale $J-J'$ devient
$$J_1-J' = \int_0^{x-xy} v_1 ^{\alpha-1} (y-x+v_1)^{\alpha-1} dv_1.$$
Ce qui devient, en posant $v_2= x-v_1$ 
$$ J_1-J' = \int_{xy}^x (x-v_2)^{\alpha-1}
 (y-v_2)^{\alpha-1} dv_2$$
 ou encore  avec $v_3=\frac{v_2}{xy} $
 $$ J_1-J' = x^\alpha y^\alpha \int_1^{1/y}(1-xv_3)^{\alpha -1}
 (1-yv_3)^{\alpha -1} dv_3;$$
 En posant pour finir 
 $ z=\frac{1}{v_3}$ on obtient le rsultat attendu, c'est  dire 
 $$ J-J' = x^\alpha y^\alpha \int_{1/y}^1 
 \frac{(z-x)^{\alpha-1} (z-y)^{\alpha-1} }{z^{2\alpha}}dz.$$
 \subsection{Calcul de la trace quand $\alpha > \demi$}
Le th\'eor\`eme \ref{inverse} et la propri\'et\'e d'uniformit\'e \'enonc\'ee dans ce th\'eor\`eme permettent 
d'\'ecrire, si $\epsilon$ est un r\'eel qui tend vers z\'ero,
\begin{eqnarray*}
\Tr \left(  T_{N}(\vert 1-\chi \vert ^{2\alpha} f_{1}) \right)^{-1}
&= & \sum _{v=N \epsilon} ^{N-N \epsilon}
 \frac{N^{2\alpha-1}}{(2\alpha-1)\Gamma ^2(\alpha) }
( \frac{k}{N}-1) ^{2\alpha-1} ( \frac{k}{N})  ^{2\alpha-1} + o(N^{2\alpha})\\
&=&  \frac{N^{2\alpha}}{(2\alpha-1)\Gamma ^2(\alpha) }
\int_{0} ^1 x^{2\alpha-1} (x-1)^{2\alpha-1} dx + o(N^{2\alpha})\\
&=&  \frac{N^{2\alpha}}{(2\alpha-1)\Gamma ^2(\alpha) }
B(2\alpha,2\alpha)  + o(N^{2\alpha})
\end{eqnarray*}

\section{Dmonstration du thorme \ref{valeurpropre}}
\subsection {Cas o l'exposant $\alpha$ est suprieur  $ \frac{1}{2}$}
Le rel $\alpha$ tant fix
nous introduisons la notation $ \varphi_\alpha = \vert 1-\chi\vert ^{2 \alpha}f_1$. 
D\'emontrons d'abord le lemme 
\begin{lemme} \label {petit}
Posons $ \varphi_\alpha = \vert 1-\chi\vert ^{2 \alpha}f_1$ avec $\alpha > \demi$. Alors  
si $\frac{k}{N}$ tend vers $0$, et $0\le l \le N$ on a  
$$\left( T_N (\varphi_\alpha) \right)^{-1}_{k+1,l+1} = o(N^{2\alpha-1})$$
$$\left( T_N (\varphi_\alpha) \right)^{-1}_{N-k+1,l+1} = o(N^{2\alpha-1}).$$
De m\^eme si $\frac{l}{N}$ tend vers $0$, et $0\le k \le N$ on a
$$\left( T_N (\varphi_\alpha) \right)^{-1}_{k+1,l+1} = o(N^{2\alpha-1})$$
$$\left( T_N (\varphi_\alpha) \right)^{-1}_{k+1,N-l+1} = o(N^{2\alpha-1}),$$
les quatre formules ci-dessus tant uniformes sur $[0,[N\epsilon]]\times [0,N]$
ou sur $[0,N] \times [0, [N\epsilon]]$ pour $\epsilon$ assez petit.
\end{lemme}
\begin{preuve} {} 
Avec les m\^emes notations que prcdemment nous avons toujours 
$$ 
(T_{N}(\varphi_\alpha))_{k+1,l+1}^{-1}=\sum_{u=0}^{\min (k,l)} \bar \gamma^\alpha_{k-u}
 \gamma^\alpha_{l-u}
- \sum_{u=N-\max(k,l)} ^{N-\max(k,l)+\min(k,l)} \bar \gamma^\alpha_{u}
 \gamma^\alpha_{u-\max(k,l)+\min(k,l)}.$$
En utilisant les points 1,2,3, de l'hypothse de rcurrence on obtient facilement 
que si $\frac{k}{N}$ ou $\frac{l}{N}$ dans $[0,N \epsilon]$ avec 
$\epsilon \rightarrow 0$ les deux sommes intervenant de l'galit ci-dessus 
sont d'ordre $o(N^{2\alpha-1})$.\\
De m\^{e}me si $\frac{k}{N} =x>0$ et $\frac{l}{N} =1-\epsilon$ on a, en rutilisant les 
calculs de la dmonstration du thorme (\ref{inverse}) :
\begin{eqnarray*}
 \sum_{u=0}^{k}&&\bar \gamma^\alpha_{k-u}  \gamma^\alpha_{l-u} =\\
&=& \frac{N^{2\alpha -1}} {\Gamma^2 f_1(1)}
  \int_0^x (x-t)^{\alpha-1} (1-x+t)^\alpha (1-t)^{\alpha-1} t^{\alpha} dt
  +o(N^{2\alpha-1})
  \end{eqnarray*}
  et
  \begin{eqnarray*}
  \sum_{u=N-l} ^{N-l+k} && \gamma^\alpha_{u}
 \bar \gamma^\alpha_{u-l+k}=\\
 &=& \frac{N^{2\alpha -1}} {\Gamma^2 f_1(1)}
 \int_0^x (1+t-x)^{\alpha-1} (x-t)^\alpha t^{\alpha-1}
 (1-t)^\alpha dt+o(N^{2\alpha-1}).
  \end{eqnarray*}
  Les m\^{e}mes changements de variables que ceux effectus dans la dmonstration du thorme \ref{inverse} donnent alors 
  $$(T_{N}(\varphi_\alpha))_{k+1,l+1}^{-1} =o(N^{2\alpha -1}).$$
  \end{preuve}
  D\'emontrons maintenant le th\'eor\`eme proprement dit dans le cas $\alpha$ sup\'erieur \`a $\demi$.
Nous avons bien videmment, pour tout entier naturel $s$  et pour tout entier $k , 0 \le k\le N$ 
\begin{eqnarray*}
\left( \left( T_N (\varphi_\alpha) \right)^{-1}\right) ^s _{k+1, k+1} &=&
\sum _{h_{s-1}=0} ^N \left( T_N (\varphi_\alpha) \right)^{-1}_{k+1,h_{s-1}+1} \left( 
 \sum _{h_{s-2}=0} ^N \left( T_N (\varphi_\alpha) \right)^{-1}_{h_{s-1}+1,h_{s-2}+1} \right.\\
&& \left. \cdots \sum _{h_{1}=0} ^N \left( T_N (\varphi_\alpha) \right)^{-1}_{h_{1}+1,k+1}\right)
 \end{eqnarray*}
Gr\^ace au thorme \ref{inverse} et au lemme (\ref{petit}) nous pouvons crire, si N suffisamment grand
 
\begin{eqnarray*}
\left( \left( T_N (\varphi_\alpha) \right)^{-1}\right) ^s _{k+1, k+1} &=&
\frac{N^{2s\alpha -s}}{\Gamma(\alpha) ^{2s} {f_1(1)}^s} \\
&\times &\sum _{h_{s-1}=0} ^N G_\alpha (k/N,h_{s-1}/N)\\
&& \times \left( 
\sum _{h_{s-2}=0} ^N G_\alpha (h_{s-1}/N,h_{s-2}/N)
\cdots  \sum _{h_{1}=0} ^N G_\alpha (h_1/N,k/N)\right) +o(N^{2s\alpha -s}).
 \end{eqnarray*}
En appliquant $s$ fois la formule d'Euler Mac-Laurin il vient 
\begin{eqnarray*}
&&\left( \left( T_N (\varphi_\alpha) \right)^{-1}\right) ^s _{k+1, k+1} =
\frac{N^{2s\alpha }}{\Gamma(\alpha) ^{2s} {f_1(1)}^s} \\
&&\times  \int_0^1 \cdots \int_0^1 G_\alpha ( k/N, t_{s-1}) 
G_\alpha (  t_{s-1},t_{s-2}) \cdots G_\alpha ( t_1, k/N) dt_{n-1} dt_{n-2} \cdots dt_1 
+o(N^{2s\alpha}) .
\end{eqnarray*}
C'est  dire que, en utilisant les notations de l'introduction nous pouvons crire 
que 
$$ \left( \left( T_N (\varphi_\alpha) \right)^{-1}\right) ^s _{k+1, k+1} =
\frac{N^{2s\alpha }}{\Gamma(\alpha) ^{2s} {f_1(1)}^s} \ast ^s G_{\alpha}(k/N,k/N) + o(N^{2s\alpha})
$$
ou encore 
$$\trace\left( \left( \left( T_N (\varphi_\alpha) \right)^{-1}\right) \right)^s  =\int_0^1 
\frac{N^{2s\alpha +1}}{\Gamma(\alpha) ^{2s} {f_1(1)}^s} \ast ^s G_{\alpha}(t,t) dt + o(N^{2s\alpha})
$$
La fonction $G_\alpha$ tant continue sur $[0,1] \times [0,1]$ on peut passer  la limite sur $s$
et crire 
\begin{eqnarray*}
& &\lim_{s\rightarrow + \infty}
\left( \trace \left( \left( T_N (\varphi_\alpha) \right)^{-1}\right) ^s \right) ^{1/s} =\\
&= & 
\frac{N^{2\alpha}}{\Gamma ^{2} f_1(1)} 
 \lim_{s\rightarrow + \infty}\left( \int_0^1\ast ^s G_{\alpha}(t,t) dt \right)^{1/s} 
+ o(N^{2\alpha}).
\end{eqnarray*}
\subsection {Cas o l'exposant $\alpha \in ]0, \frac{1}{2}[$}
Remarquons que comme pour $\alpha<\demi$ nous avons le lemme suivant 
qui correspond au lemme (\ref{petit}).
\begin{lemme} \label {petitpetit}
Posons $ \varphi_\alpha = \vert 1-\chi\vert ^{2 \alpha}f_1$ avec $\alpha \in]0, \demi[$. Alors  
si $\frac{k}{N}$ tend vers $0$, et 
 $\frac{l}{N} \rightarrow y >0$ on a 
$$\left( T_N (\varphi_\alpha) \right)^{-1}_{k+1,l+1} = o(N^{2\alpha-1}),$$
si, d'autre part, $\frac{l}{N} \rightarrow y<1$  on a
$$\left( T_N (\varphi_\alpha) \right)^{-1}_{N-k+1,l+1} = o(N^{2\alpha-1}).$$
De m\^eme si $\frac{l}{N}$ tend vers $0$, et $\frac{k}{N} \rightarrow x>0$ on a
$$\left( T_N (\varphi_\alpha) \right)^{-1}_{k+1,l+1} = o(N^{2\alpha-1})$$
si, d'autre part, $\frac{k}{N} \rightarrow x<1$ on a
$$\left( T_N (\varphi_\alpha) \right)^{-1}_{k+1,N-l+1} = o(N^{2\alpha-1}),$$
les quatre formules ci-dessus tant uniformes sur $[0,[N\epsilon]]\times [0,N]$
ou sur $[0,N] \times [0, [N\epsilon]]$ pour $\epsilon$ assez petit.
\end{lemme}

Pour pouvoir reprendre la d\'emonstration du cas $\alpha> \frac{1}{2}$ nous devons \'etablir les trois lemmes suivants.
\begin{lemme}\label{convergence1}
Si $0<\alpha<\frac{1}{2}$ et $0< x\not= y <1$ on a 
$$ G_{\alpha} (x,y) \le  
 \frac{\Gamma(1-2\alpha) \Gamma (\alpha)}{\Gamma(1-\alpha)}
 \left(\max (x,y) -\min (x,y)\right)^{2\alpha-1}.$$
\end{lemme}
\begin{preuve}{}
Pour la d\'emonstration de ce lemme nous supposerons $0<x<y<1$. Il est alors clair que 
$$\int_{y}^1 \frac{(t-x)^{\alpha-1} (t-y)^{\alpha-1} }{t^{2\alpha}} dt \le \frac{1} {y^{2\alpha}}
\int_{y}^1 (t-x)^{\alpha-1} (t-y)^{\alpha-1}  dt.$$
Nous allons nous concentrer sur l'int\'egrale 
 $\int_{y}^1 (t-x)^{\alpha-1} (t-y)^{\alpha-1}  dt.$
 En utilisant des changements de variables successifs nous obtenons 
\begin{eqnarray*}
\int_{y}^1 (t-x)^{\alpha-1} (t-y)^{\alpha-1} dt &=& \int_{y-x}^{1-x} h ^{\alpha-1} (\frac{h}{y-x} -1)^{\alpha-1} dh
 \, (y-x)^{\alpha-1} \\
&=& \int_{1}^{(1-x)/(y-x)} u^{\alpha-1} (u-1)^{\alpha-1} du \,(y-x) ^{2\alpha-1} \\
&=& \int_{(y-x)/(1-x)}^1 v^{-2\alpha} (1-v) ^{\alpha-1} dv \, (y-x) ^{2\alpha-1}
\end{eqnarray*}
Nous savons d'autre part que 
$$ \int_{0}^1 v^{-2\alpha} (1-v) ^{\alpha-1} dv 
= \frac{\Gamma(1-2\alpha) \Gamma (\alpha)}{\Gamma(1-\alpha)}.$$
Ce qui permet finalement de conclure 
$$ G_{\alpha}(x,y) \le \frac{x^\alpha}{y^\alpha} 
\frac{\Gamma(1-2\alpha) \Gamma (\alpha)}{\Gamma(1-\alpha)}  (y-x) ^{2\alpha-1},$$
et finalement, puisque $\alpha$ est positif 
$$ G_{\alpha}(x,y) \le  
\frac{\Gamma(1-2\alpha) \Gamma (\alpha)}{\Gamma(1-\alpha)}  (y-x) ^{2\alpha-1}.$$
\end{preuve}
\begin {lemme} \label{convergence2}
En posant $t_{\alpha}(x,y) = \left( \max (x,y) - \min (x,y) \right)^{2\alpha-1} $ 
on a, si $0< y \not= z <1$ 
$$ \int_{0}^1 t_{\alpha}(y,x) t_{\alpha}(x,z) dt \le K_{\alpha} t_{\alpha}(y,z). $$
avec $K_{\alpha}= \left (\frac{1}{2\alpha} +\frac{\Gamma^2(2\alpha)}{\Gamma(4\alpha) } 
+\frac{\Gamma(1-2\alpha) \Gamma (\alpha)}{\Gamma(1-\alpha)}\right)$.
\end{lemme}
\begin{preuve}{}
Supposons $0<y<z<1$. Posons 
\begin{eqnarray*}
\int_{0}^1 t_{\alpha}(y,x) t_{\alpha}(x,z) dx &=& \int_{0}^y t_{\alpha}(y,x) t_{\alpha}(x,z) dx +\\
&+&\int_{y}^z t_{\alpha}(y,x) t_{\alpha}(x,z) dx + \int_{z}^1 t_{\alpha}(y,x) t_{\alpha}(x,z) dx.
\end{eqnarray*}
 Avec des changements de variables du m\^eme style que ceux 
utilis\'es pr\'ec\'edemment  nous obtenons les \'egalit\'es
\begin{eqnarray*}
\int_{0}^y t_{\alpha}(y,x) t_{\alpha}(x,z) dx &=&  \int _{z-y}^z u^{2\alpha-1} (y-z+u)^{2\alpha-1} du\\
&=& \int_{z-y}^z u^{2\alpha-1} \left( \frac{u}{z-y} -1\right)^{2\alpha-1} du (z-y) ^{2\alpha-1}\\
&=&  \int _{1}^{z/(z-y)} v^{2\alpha-1} (v-1)^{2\alpha-1} dv (z-y)^{4\alpha-1}. 
\end{eqnarray*}
On peut remarquer que 
$$
\int _{1}^{z/(z-y)} v^{2\alpha-1} (v-1)^{2\alpha-1} dv \le 
\int _{1}^{z/(z-y)}  (v-1)^{2\alpha-1} dv.
$$
Et puisque 
$$
\int _{1}^{z/(z-y)}  (v-1)^{2\alpha-1} dv = \frac{1}{2\alpha} \frac{y^{2\alpha} }{(z-y)^{2\alpha}},
$$
nous pouvons finalement conclure 
$$
\int_{0}^y t_{\alpha}(y,x) t_{\alpha}(x,z) dx  \le \frac{1}{2\alpha} t_{\alpha}(y,z).
$$
 Nous obtenons de m\^{e}me
\begin{eqnarray*}
\int_{y}^z t_{\alpha}(y,x) t_{\alpha}(x,z) dx &=& \int _{0}^{z-y} (z-y-u)^{2\alpha-1} u^{2\alpha-1} du \\
&=& \int_{0}^{z-y}  \left( 1 - \frac{u}{z-y} \right) ^{2\alpha-1} u^{2\alpha-1} du (z-y)^{2\alpha-1}\\
&=&  (z-y)^{4\alpha-1}\int _{0}^1
(1-v)^{2\alpha-1} v^{2\alpha-1} dv = \frac{\Gamma^2(2\alpha)}{\Gamma(4\alpha) } (z-y)^{4\alpha-1}.
\end{eqnarray*}
Puisque $\alpha>0$ nous pouvons conclure 
$$\int_{y}^z t_{\alpha}(y,x) t_{\alpha}(x,z) dx \le \frac{\Gamma^2(2\alpha)}{\Gamma(4\alpha) } 
t_{\alpha}(y,z).$$
Enfin l'int\'egrale $ \int_{z}^1 t_{\alpha}(y,x) t_{\alpha}(x,z) dx$ a d\'eja \'et\'e \'etudi\'ee pr\'ec\'edemment.
\end{preuve}
Le lemme suivant donne une  expression asymptotique de $\left( T_{N} \left( \vert 1- \chi \vert ^{2\alpha} f_{1}\right)^{-1}\right)_{k+1,l+1}$
pour   $\alpha\in ]0, \frac{1}{2}[$ qui est plus compliqu\'ee que la pr\'ec\'edente mais qui a l'int\'er\^{e}t d'\^{e}tre d\'efinie pour $x=y$ 
sous l'hypoth\`ese $\alpha<\demi$.
\begin{lemme}\label{convergence3}
Soient $\alpha \in ]0, \frac{1}{2}[$ et $x,y$ deux r\'eels tels que $0 < x ,y <1$. On a, en posant $f= \vert 1- \chi \vert ^{2\alpha} f_{1}$, et $\mu_{\alpha}(x,y,t)=(\max(x,y)-\min(x,y)+t)$ 
\begin{eqnarray*} 
\left( T_{N}(\varphi_\alpha) \right)^{-1}_{[Nx]+1,[Ny]+1} &=& \widehat{\frac{1}{\varphi_\alpha}} \Bigl(\max ([Nx],[Ny])- \min ([Nx],[Ny]) \Bigr)+\\
&+& \frac{N^{2\alpha-1}}{\Gamma^2(\alpha) f_{1}(1)}h_{\alpha}(x,y)+o(N^{2\alpha-1}),
\end{eqnarray*}
avec :
\begin{eqnarray*}
h_{\alpha}(x,y)&=& \int_{0}^{\min (x,y)} t^{\alpha-1} (\mu_{\alpha}(x,y,t))^{\alpha-1} \Bigl( 2(1-t)^\alpha
 (1-\mu_{\alpha}(x,y,t))^\alpha -(1-t)^\alpha -(1-\mu_{\alpha})^\alpha\Bigr) dt \\
&-& \int_{1-\min(x,y)}^1t^{\alpha-1} (1-t)^\alpha\Bigl(2 t- \mu_{\alpha}(x,y,t)\Bigr)^{\alpha-1}
 \Bigl(1- \mu_{\alpha}(x,y,t)\Bigr)^\alpha dt\\
&-& \int_{\min(x,y)}^{+ \infty} t^{\alpha-1} \left(\mu_{\alpha}(x,y,t)\right)^{\alpha-1}dt
\end{eqnarray*}
\end{lemme}
\begin{preuve}{}
Dans cette d\'emonstration nous supposerons que $x<y$. Nous allons de nouveau utiliser la formule de Gohberg-Semencul, mais nous allons regrouper les
termes dans un ordre diff\'erent de ce qui a \'et\'e fait dans la d\'emonstration du th\'eor\`eme 
\ref{inverse}. On a, toujours avec les notations introduites au cours de l'article, et en posant,
 pour all\'eger les notations, $ [Nx] =k$, $[Ny] =l$, 
$$ 
\left(T_{N}\right)^{-1}_{k+1,l+1} = 
\left( \bar\gamma_{0}^\alpha \gamma_{l-k}^{(\alpha)} + \cdots + \bar\gamma_{k}^{(\alpha)} 
\gamma_{l}^{(\alpha)}\right)
-\left ( \bar \gamma_{N-k}^{(\alpha)} \gamma_{N-l}^{(\alpha)}+ \cdots \bar \gamma_{N}^{(\alpha)}
 \gamma_{N-k+l}^{(\alpha)}
\right)$$
Nous remarquons qu'en normalisant par $\vert \beta_{0}^{(\alpha)} \vert ^2 $  il vient 
\begin{eqnarray*}
\sum_{u=0}^k \bar \gamma_{u}^\alpha \gamma_{l-k+u}^\alpha &=&
\sum_{u=0}^k \bar \beta_{u}^{(\alpha)} \beta_{l-k+u}^{(\alpha)}+\\
&+& \sum_{u=0}^k \bar \beta_{u}^{(\alpha)} ( \gamma^{(\alpha)} _{l-k+u} - \beta_{l-k+u}^{(\alpha)})+\\
&+& \sum_{u=0}^k \overline{(\gamma_{u}^{(\alpha)}-\beta_{u}^{(\alpha)})} \gamma^{(\alpha)} _{l-k+u} +\\
&+&  \sum_{u=0}^k \overline{(\gamma_{u}^{(\alpha)}-\beta_{u}^{(\alpha)})}
 ( \gamma^{(\alpha)} _{l-k+u} - \beta_{l-k+u}^{(\alpha)}).
 \end{eqnarray*}
Nous obtenons dans un premier temps 
\begin{eqnarray*}
\sum_{u=0}^k \bar \beta_{u}^{(\alpha)} \beta_{l-k+u}^{(\alpha)}
&=& 
\sum_{u=0}^{+\infty} \bar \beta_{u}^{(\alpha)} \beta_{l-k+u}^{(\alpha)} 
- \sum_{u=k+1}^{+\infty} \bar \beta_{u}^{(\alpha)} \beta_{l-k+u}^{(\alpha)} \\
&=&  \widehat{\frac{1}{\varphi_\alpha}} \left(l-k \right) -\frac{N^{2\alpha-1}}{\Gamma^2(\alpha) f_{1}(1)} 
I_{1,\alpha} (x,y) + o(N^{2\alpha-1}),
\end{eqnarray*}
puis 
\begin{eqnarray*}
\sum_{u=0}^k \bar \beta_{u}^{(\alpha)} ( \gamma^{(\alpha)} _{l-k+u} - \beta_{l-k+u}^{(\alpha)})
&=&   \frac{N^{2\alpha-1}}{\Gamma^2(\alpha) f_{1}(1)}I_{2\alpha}(x,y) + o(N^{2\alpha-1}) \\
 \sum_{u=0}^k \overline{(\gamma_{u}^{(\alpha)}-\beta_{u}^{(\alpha)})} \gamma^{(\alpha)} _{l-k+u}
 &=&  \frac{N^{2\alpha-1}}{\Gamma^2(\alpha) f_{1}(1)}I_{3\alpha} (x,y)+ o(N^{2\alpha-1}) \\
 \sum_{u=0}^k \overline{(\gamma_{u}^{(\alpha)}-\beta_{u}^{(\alpha)})}
 ( \gamma^{(\alpha)} _{l-k+u} - \beta_{l-k+u}^{(\alpha)})&= &
 \frac{N^{2\alpha-1}}{\Gamma^2(\alpha) f_{1}(1)} I_{4\alpha} (x,y)+ o(N^{2\alpha-1}) 
\end{eqnarray*}
avec 
\begin{eqnarray*}
I_{1\alpha} (x,y)
&=& \int _{x}^{+ \infty} t^{\alpha-1} (y-x+t)^{\alpha-1} dt \\
I_{2\alpha} (x,y)
&=&
\int _{0}^{x} t^{\alpha-1} (y-x+t) ^{\alpha-1} \left( (1-y+x-t)^\alpha -1\right)dt \\
 I_{3\alpha} (x,y)&=& 
 \int_{0}^x
 (y-x+t)^{\alpha-1} (1-y+x-t)^{\alpha} t^{\alpha-1} \left( (1-t)^\alpha-1 \right) dt \\
 I_{4\alpha} (x,y) &= &
\int_{0}^x  
 t^{\alpha-1}( y-x+t)^{\alpha-1} \left( (1-t)^\alpha -1\right) \left( (1-y+x-t)^\alpha-1\right)dt.
\end{eqnarray*}
Nous avons d'autre part \'etabli dans un pr\'ec\'edent travail (\cite{RS10} ) que 
$$\sum_{u=0}^k  \gamma^{(\alpha)} _{N-k+u}\bar \gamma^{(\alpha)} _{N-l+u} 
=\frac{N^{2\alpha-1}}{\Gamma^2(\alpha) f_{1}(1)}
 I_{5,\alpha} (x,y)
+ o(N^{2\alpha-1}),$$ avec 
$$ I_{5,\alpha} (x,y) =  \int_{1-x}^1 t^{\alpha-1} (1-t)^{\alpha} (t-y+x)^{\alpha-1} (1-t+y-x)^\alpha dt.$$
Tout ceci donne finalement comme annonc\'e 
$$ 
\left( T_{N}(f) \right)^{-1}_{k+1,l+1} = \widehat{\frac{1}{f}} \left(l-k) \right)
+ \frac{N^{2\alpha-1}}{\Gamma^2(\alpha) f_{1}(1)}h_{\alpha}(x,y)+o(N^{2\alpha-1}),
$$
avec 
$$ h_{\alpha}(x,y) = I_{2,\alpha}(x,y) +I_{3,\alpha}(x,y) + I_{4,\alpha}(x,y) -I_{1,\alpha}(x,y) -I_{5,\alpha}(x,y) .$$
En rempla\c cant les int\'egrales par leur valeur  on obtient l'expression de 
$h_{\alpha}(x,y)$ \'enonc\'ee dans le lemme.
\end{preuve}

Nous allons maintenant passer \`a la d\'emonstration du th\'eor\`eme proprement dit.\\
Nous pouvons bien sr \'ecrire 
\begin{align} 
\trace \left ( \left( (T_{N}(\varphi_\alpha))^{-1}\right)^{s} \right) 
&= 
\sum _{0 \le k_{1}\le N }\left( (T_{N}(\varphi_\alpha))^{-1}\right)_{k_{1}+1, k_{2}+1}
\sum _{0 \le k_{2}\le N }\left( (T_{N}(\varphi_\alpha))^{-1}\right)_{k_{2}+1, k_{3}+1}\\ \nonumber 
&\cdots
 \sum _{0 \le k_{s}\le N }\left( (T_{N}(\varphi_\alpha))^{-1}\right)_{k_{s-1}+1, k_{s}+1}
\left( (T_{N}(\varphi_\alpha))^{-1}\right)_{k_{s}+1, k_{1}+1}.\nonumber
\end{align}
En utilisant le  lemme (\ref{convergence3}) et le lemme (\ref{petitpetit}) nous obtenons,si $\epsilon$ d\'esigne un r\'eel positif,
 \begin{eqnarray*} 
 & & \sum _{0 \le k_{s}\le N }\left( (T_{N}(\varphi_\alpha))^{-1}\right)_{k_{s-1}+1, k_{s}+1} 
\left( (T_{N}(\varphi_\alpha))^{-1}\right)_{k_{s}+1, k_{1}+1} = \\
&=&\frac{N^{4\alpha-1}}{\Gamma^4(\alpha) f^2_{1}(1)} \left( 
\int_{t_{s}\in\left( [0,1] \setminus \Delta_{\epsilon,1,s}\right)}  
G_{\alpha} (k_{s-1}/N, t_{s}) G_{\alpha} (t_{s},k_{1}/N) \, dt_{s}Ê+\right . \\
&+&\left . \int_{t_{s}\in \Delta_{\epsilon,1,s}} 
h_{\alpha} (k_{s-1}/N, t_{s}) h_{\alpha} (t_{s},k_{1}/N) \,dt_{s} +\right.\\
&+& \left.
 \sum _{k_{s}\in D_{\epsilon,1,s}}
\widehat{\frac{1}{\varphi_\alpha}} \left( \max(k_{s-1},k_{s}) - \min  (k_{s-1},k_{s})\right) 
\widehat{\frac{1}{\varphi_\alpha}} \left( \max(k_{1},k_{s}) - \min  (k_{1},k_{s})\right) \right) +\\
&+& o(N^{4\alpha-1})
\end{eqnarray*}
o\`u 
$ \Delta_{\epsilon,1,s}=   [\frac{k_{1}}{N}-\epsilon, \frac{k_{1}}{N}+\epsilon] 
\bigcup [\frac{k_{s-1}}{N}-\epsilon, \frac{k_{s-1}}{N}+\epsilon]$\\
et $D_{[N\epsilon],1,s-1}= \{ k\in \mathbb N / k_{s-1} -[N \epsilon ] \le k
\le  k_{s-1}  +[N \epsilon], k_{1} -[N \epsilon] \le k\le k_{1}   +[ N\epsilon] \},$ o\`u $[t]$ d\'esigne la partie enti\`ere d'un r\'eel $t$. \\
 En remarquant que 
si $y,z \in [0,1]$ la fonction $x \rightarrow G_{\alpha} (y,x)G_{\alpha} (x,z)$ est int\'egrable sur $[0,1]$  
(lemme (\ref{convergence2})), on obtient, en faisant tendre $\epsilon$ vers z\'ero 
 \begin{eqnarray*} 
 & & \sum _{0 \le k_{s}\le N }\left( (T_{N}(\varphi_\alpha))^{-1}\right)_{k_{s-1}+1, k_{s}+1} 
\left( (T_{N}(\varphi_\alpha))^{-1}\right)_{k_{s}+1, k_{1}+1} = \\
&=& \frac{N^{4\alpha-1}}{\Gamma^4(\alpha) f^2_{1}(1)} 
\int_0^1 
G_{\alpha} (k_{s-1}/N, t_{s}) G_{\alpha} (t_{s},k_{1}/N) \, dt_{s}Ê+
 o(N^{4\alpha-1})
\end{eqnarray*}
si $\alpha>\frac{1}{4}$ et, dans le cas contraire 
 \begin{eqnarray*} 
 & & \sum _{0 \le k_{s}\le N }\left( (T_{N}(\varphi_\alpha))^{-1}\right)_{k_{s-1}+1, k_{s}+1} 
\left( (T_{N}(\varphi_\alpha))^{-1}\right)_{k_{s}+1, k_{1}+1} = \\
&=&\frac{N^{4\alpha-1}}{\Gamma^4(\alpha) f^2_{1}(1)} 
\int_0^1
G_{\alpha} (k_{s-1}/N, t_{s}) G_{\alpha} (t_{s},k_{1}/N) \, dt_{s}Ê+ \\
&+& 
 \sum _{k_{s}\in D_{n_{0},1,s-1}}
\widehat{\frac{1}{\varphi_\alpha}} \left( \max(k_{s-1},k_{s}) - \min  (k_{s-1},k_{s})\right) 
\widehat{\frac{1}{\varphi_\alpha}} \left( \max(k_{1},k_{s}) - \min  (k_{1},k_{s})\right)  +\\
&+& o(N^{4\alpha-1})
\end{eqnarray*}
si $n_{0}$ est un entier tel que si $ \vert m \vert >n_{0}$ on peut remplacer $\hat \varphi_{\alpha}(m)$ par 
son asymptotique. Puisque $\alpha$ est un r\'eel positif on obtient, pour un entier $i$ suffisamment grand en r\'ep\'etant \`a chaque \'etape le m\^eme d\'ecoupage, 
   \begin{eqnarray*} 
\trace \left ( \left( (T_{N}(\varphi_\alpha))^{-1}\right)^{s} \right) 
&= &
\sum _{0 \le k_{1}\le N }\left( (T_{N}(\varphi_\alpha))^{-1}\right)_{k_{1}+1, k_{2}+1}
\sum _{0 \le k_{2}\le N }\left( (T_{N}(\varphi_\alpha))^{-1}\right)_{k_{2}+1, k_{3}+1}\\ 
&& \cdots \sum _{0 \le k_{i-1}\le N }\left( (T_{N}(\varphi_\alpha))^{-1}\right)_{k_{i-1}+1, k_{i}}\\ 
&& \frac{(N^{2\alpha(i+2)-1} }{\left(\Gamma^2(\alpha) f_{1}(1)\right)^{i+1}} 
\int_0 ^1 G_{\alpha} (k_{i-1}/N, t_{i})  dt_{i}\cdots \\
&& \cdots\frac{N^{4\alpha-1}}{\Gamma^4(\alpha) f^2_{1}(1)} 
\int_0 ^1 
G_{\alpha} (t_{s-1}, t_{s}) G_{\alpha} (t_{s},k_{1}/N) dt_{s} 
\left( 1+o(1)\right),
\end{eqnarray*}
En r\'ep\'etant cette m\'ethode on peut de m\^eme obtenir pour tout entier $i$ $1\le i \le s-1$,
le passage de l'expression ci-dessus \`a l'expression 
  $$
\trace \left ( \left( (T_{N}(\varphi_\alpha))^{-1}\right)^{s} \right) 
=  \frac{N^{2 s \alpha-1}}{\Gamma^{2s}(\alpha) f^{2s}_{1}(1)} 
\int_0 ^1 *^s G_{\alpha} (t_{1}, t_{1}) dt_{1} 
\left( 1+o(1)\right).
$$

 \section{Dmonstration du corollaire \ref{encadrements}}
 \subsection{Cas o l'exposant $\alpha$ est dans $]0, \demi[$}
 Nous allons d'abord encadrer $\Lambda_{\alpha,N}$ la plus grande valeur propre de
 $\left(T_N (\varphi_\alpha)\right)^{-1}$. Pour cela nous allons encadrer 
 $\trace \left(\left(\left(T_N (\varphi_\alpha)\right)^{-1}\right)^s\right) $.\\
 \textbf{Majoration de $\Lambda_{\alpha,N}
 $ avec  $\alpha \in ]0, \demi [$}\\
 Une majoration est immdiatement fournie par les lemmes 
 \ref{convergence1} \ref{convergence2}. Cette majoration est 
 $$ \trace \left(\left(\left(T_N (\varphi_\alpha)\right)^{-1}\right)^s\right)
 \le N^{2\alpha s}
\left(  \frac{\Gamma (1-2\alpha) \Gamma (\alpha)} {\Gamma (1- \alpha)}\right)^s
 K_\alpha ^{s}\left( \frac{1}{ \Gamma ^2 (\alpha) f^2_1 (1)}
 \right)^s\int_0^1 \int_0 ^1 t_\alpha (t_1,t_2)  dt_1 dt_2 .$$ 
 Soit encore 
 $$ \trace \left(\left(\left(T_N (\varphi_\alpha)\right)^{-1}\right)^s\right)
 \le N^{2\alpha s}
  \frac{\Gamma (1-2\alpha) \Gamma (\alpha)} {\Gamma (1- \alpha)}
 K_\alpha ^{s}\left( \frac{1}{ \Gamma ^2 (\alpha) f^2_1 (1)}
 \right)^s \frac {1}{\alpha(2\alpha+1)} .$$ 
 Puisque nous avons 
 $$ \Lambda_{\alpha,N} = \lim_{s \rightarrow + \infty} 
\Bigl(\trace \left(\left(\left(T_N (\varphi_\alpha)\right)^{-1}\right)^s\right)\Bigr)^{1/s},$$
nous obtenons la majoration
$$ \Lambda_{\alpha,N} \le N^{2\alpha} K_\alpha  \frac{1}{ \Gamma ^2 (\alpha) f^2_1 (1)} 
\frac{\Gamma (1-2\alpha) \Gamma (\alpha)} {\Gamma (1- \alpha)}$$
qui s'crit aussi, en remplaant $K_\alpha$ par sa valeur, 
 $$ 3 \Lambda_{\alpha,N} \le N^{2\alpha} 
  \left( \frac{1}{2\alpha}, \frac{\Gamma ^2(2\alpha)}{\Gamma (4 \alpha)}
  , \frac{\Gamma (1-2\alpha) \Gamma (\alpha)}{\Gamma (1-\alpha)}\right)
   \frac{1}{ \Gamma  (\alpha) f^2_1 (1)} \frac{\Gamma (1-2\alpha)}{\Gamma(1-\alpha)}$$
 \textbf{Minoration de $\Lambda_{\alpha,N}
 $ avec  $\alpha \in ]0, \demi [$.}\\
 Supposons ici que $0<x<y<1$. Nous avons 
 $$ G_{\alpha}(x,y) \ge x^\alpha y^\alpha \int_{y}^1 (t-x)^{\alpha-1} (t-y)^{\alpha-1} dt ,$$
 soit en reprenant l'expression obtenue dans la d\'emonstration du lemme \ref{convergence1} :
 \begin{eqnarray*} 
 G_{\alpha}(x,y) &\ge& \frac{\Gamma (1-2\alpha) \Gamma(\alpha) }{\Gamma(1-\alpha)} x^\alpha y^\alpha (y-x)^{2\alpha-1}\\
 &\ge &\frac{\Gamma (1-2\alpha) \Gamma(\alpha) }{\Gamma(1-\alpha)} x^{2\alpha} (y-x)^{2\alpha-1}\\
 & \ge& \frac{\Gamma (1-2\alpha) \Gamma(\alpha) }{\Gamma(1-\alpha)} x^{2\alpha} y^{2\alpha-1}\\
 & \ge& \frac{\Gamma (1-2\alpha) \Gamma(\alpha) }{\Gamma(1-\alpha)} x^{2\alpha} y^{2\alpha}.
 \end{eqnarray*}
Ce qui donne 
 $$ \trace \left(\left(\left(T_N (\varphi_\alpha)\right)^{-1}\right)^s\right )
  \ge  \frac{ N^{2 s \alpha } }
  {(\Gamma ^2 (\alpha) f_1(1) ) ^s} \left(\frac{\Gamma (1-2\alpha) \Gamma(\alpha) }{\Gamma(1-\alpha)} \right)^s
\left(\int_{0}^1 t^{4\alpha} dt\right)^s.$$
Ce qui nous donne finalement 
$$
\Lambda_{\alpha,N} \ge \frac{N^{2\alpha}}{f_{1}(1)} \frac{\Gamma(1-2\alpha) }{\Gamma(1-\alpha) \Gamma(\alpha) (4\alpha+1)}.$$
  \subsection{Cas o l'exposant $\alpha$ est dans $]\demi,1[$}
 L aussi nous allons chercher  encadrer $\Lambda_{\alpha,N}$.\\
  \textbf{Majoration de $\Lambda_{\alpha,N}
 $ avec  $\alpha \in ]\demi , 1[$.}\\
 En utilisant la formule de trace obtenue plus haut nous obtenons facilement 
 $$\Lambda_{\alpha,N} \le N^{2\alpha} \frac{B(2\alpha,2\alpha)}{\Gamma^2(\alpha) (2\alpha-1)}+o(N^{2\alpha}).$$
  \textbf{Minoration de $\Lambda_{\alpha,N}$ avec $\alpha \in ]\demi, 1[$.}\\
  Pour simplifier les notations nous supposerons dans ce paragraphe 
que $0<x<y<1$.
 Remarquons que nous avons 
 \begin{eqnarray*}
  G_\alpha (x,y) &=& x^\alpha y^\alpha 
  \int_y^1 \frac{ (t-x)^{\alpha-1}(t-y)^{\alpha-1}}{t^{2\alpha} }dt \\
  & \ge &  x^\alpha y^\alpha 
  \int_y^1 (t-x)^{\alpha-1} dt \\
  &\ge& x^\alpha y^\alpha  \left( \frac{(1-x)^\alpha -(y-x)^\alpha }{\alpha}\right).
  \end{eqnarray*}
  D'autre part d'aprs le thorme des accroissements finis 
  il existe un rel $c$, $ y<c<1$ tel que 
  $(1-x)^\alpha -(y-x)^\alpha = \alpha (1-y) (c-x)^{\alpha -1}$.
  Puisque $ (c-x)^{\alpha-1} > (1-x)^{\alpha-1}> (1-x)$ ($\alpha <1$)
  nous pouvons crire les ingalits 
   \begin{eqnarray*} 
  G_\alpha (x,y) &\ge & x^\alpha y^\alpha  (1-y) (c-x)^{\alpha -1}\\
  & \ge & x^\alpha y^\alpha  (1-y) (1-x)^{\alpha -1}\\
  & \ge & x^\alpha y^\alpha  (1-y) (1-x)
  \end{eqnarray*}
  En remarquant que 
  \begin{equation} \label{produit}
   \int_0^1 \cdots \int_0^1 g(x_1) f(x_2) g(x_2) f(x_3) 
  \cdots g(x_n) f(x_1) dx_n dx_{n-1} \cdots dx_{n-1} 
  = \left( \int_0^1 g(x) f(x) dx \right)^s 
  \end{equation}
  nous pouvons conclure que 
  $$ \trace \left(\left(\left(T_N (\varphi_\alpha)\right)^{-1}\right)^s\right )
  \ge  \frac{ N^{2 s \alpha } }
  {(\Gamma ^2 (\alpha) f_1(1) ) ^s}
   \left( \int_0^1 x^{2 \alpha} (1-x) ^2 dx \right)^s.$$
 Soit $$ \Lambda _{\alpha,N} \ge N^{2\alpha}
 \frac{1}{\Gamma ^2 (\alpha) f_1(1)}
 \int_0^1 x^{2 \alpha} (1-x) ^2 dx.$$
 Ce qui donne, tous calculs faits 
$$ \Lambda _{\alpha,N} \ge N^{2\alpha} \frac{2}{(2\alpha+1)(2\alpha+2)(2\alpha+3)}
\frac{1}{\Gamma^2(\alpha) f_{1}(1)}.$$
\subsection{Cas o l'exposant $\alpha$ est dans $]1,+ \infty[$}
   \textbf{Majoration de $\Lambda_{alpha,N}
 $ avec  $\alpha \in ]1 , +\infty[$.}\\
 Notre calcul s'inspire de celui
 effectu dans le cas d'un exposant entier par 
 B\"{o}ttcher et Widom dans \cite{BoW}. 
 Dans l'criture 
\begin{align} \label{majo}
\trace\left(\left( T_N(\varphi_\alpha) \right)^{-1} \right)^s &=  
\frac{N^{2\alpha s}}{ \left( \Gamma^2 (\alpha)f_1(1)\right)^s} \times \\
&\times   \int _0^1 \int _0^1G_\alpha (t_1,t_2) \int_0^1 G_\alpha (t_2,t_3)
\cdots \int _0^1G_\alpha (t_{s-1},t_{s}) \int_0^1 G_\alpha (t_s,t_1)
dt_1 \cdots dt_s \nonumber
\end{align}
nous posons le changement de variable 
$$ t_1 =\frac{1+x_1}{2} \, t_1 =\frac{1+x_2}{2} 
\cdots  t_s =\frac{1+x_s}{2}.$$
Nous obtenons alors 
\begin{align} \label{major}
\trace\left(\left( T_N(\varphi_\alpha) \right)^{-1} \right)^s &=  
\frac{N^{2\alpha s}}{ \left( \Gamma^2 (\alpha)f_1(1)\right)^s} \frac{1}{2^s} 
   \int _{-1}^1 \int _{-1}^1G_\alpha ( \frac{1+x_1}{2},\frac{1+x_2}{2})\times \nonumber  \\
&\times \int_{-1}^1 G_\alpha (\frac{1+x_2}{2},\frac{1+x_3}{2})
\cdots \int _{-1}^1G_\alpha (\frac{1+x_{s-1}}{2},\frac{1+x_s}{2}) \times  \\
&\times \int_{-1}^1 G_\alpha (\frac{1+x_s}{2},\frac{1+x_1}{2})
dx_1 \cdots dx_s \nonumber
\end{align}
Nous pouvons alors remarquer que 
\begin{eqnarray*}
G_\alpha (\frac{1+x_1}{2},\frac{1+x_2}{2} ) &=&
\left(\frac{1+x_1}{2}\right)^\alpha \left(\frac{1+x_2}{2}\right)^\alpha\\
&&\int_{\max \left((1+x_1)/2, (1+x_2)/2\right)} ^1 
\frac{\left(t-\frac{1+x_1}{2}\right)^{\alpha-1} 
\left(t-\frac{1+x_2}{2}\right)^{\alpha-1}}{t^{2\alpha}} dt\\
&=& \frac{(1+x_1)^\alpha (1+x_2)^\alpha}{4^\alpha}
 \int _{\max(x_1,x_2)}^1 \frac{\left(\frac{t'-x_1}{2}\right)^{\alpha-1}
\left( \frac{t'-x_2}{2}\right)^{\alpha-1}}{\left(\frac{1+t'}{2}\right)^{2\alpha}}
 \frac{dt'}{2}\\
 &=& \frac{2}{4^{2\alpha}} (1+x_1)^\alpha (1+x_2)^\alpha
 \int _{\max(x_1,x_2)}^1 \frac{\left(t'-x_1\right)^{\alpha-1}
\left( t'-x_2\right)^{\alpha-1}}{\left(\frac{1+t'}{2}\right)^{2\alpha}} dt'.
\end{eqnarray*}
 D'autre part la drive logarithmique de la fonction
 $ t \rightarrow \frac{(t-x_1) (t-x_2)}{\left(\frac{ 1+t}{2}\right)^2}$
 est $ \frac{(2+x_1+x_2)t-x_1-x_2-2x_1x_2}{(1+t)(t-x_1)(t-x_2)}$ qui est 
 positive sur $[\max (x_1, x_2), 1]$. Nous pouvons donc crire
 \begin{eqnarray*}
 G_\alpha (\frac{1+x_1}{2},\frac{1+x_2}{2} ) 
 &\le &\frac{2}{4^{2\alpha}} 
 (1+x_{1})^\alpha (1+x_{2})^\alpha
  \int_{\max (x_{1},x_{2})}^1 \frac{ (1-x_{1})^{\alpha-1} (1-x_{2}) ^{\alpha-1}}{\left( \frac{1+t}{2} \right)^2}  dt \\
 &\le & \frac{4}{4^{2\alpha}} (1+x_{1}) (1+x_{2}) (1-x_{1}^2)^{\alpha-1} (1-x_{2}^2) ^{\alpha-1}
 \end{eqnarray*}
 L'\'egalit\'e \ref{majo} fournit alors 
  $$
  \lim_{s \rightarrow + \infty} \left(\trace\left(\left(\left( T_N(\varphi_\alpha) \right)^{-1} \right)\right)^s\right)^{1/s}
 =  \Lambda_{\alpha,N} $$
 et la majoration 
$$  \Lambda_{\alpha,N}  \le  \int _{-1}^1 (1+x)^2 (1-x)^{2\alpha-2} dx 
 \frac{1}{ \Gamma ^2(\alpha) f_{1}(1) } \frac{2}{4^{2\alpha}}.
$$
 Soit 
 $$ \Lambda_{\alpha,N}  \le \frac{1}{f_{1}(1)} \frac{\Gamma (2\alpha+1) \Gamma (2\alpha-1) }
 { \Gamma(4\alpha) \Gamma^2 (\alpha)}.$$
 Ce qui donne une minoration de $c_{\alpha}$ \'equivalente \`a celle donn\'ee dans l'article \cite{BoW} pour le cas entier.\\
 \textbf{Minoration de $\Lambda_{\alpha,N}
 $ avec  $\alpha \in ]1 , +\infty[$.}\\
Cette fois-ci encore nous supposerons $0<x<y<1$. Nous pouvons alors \'ecrire 
\begin{eqnarray*} 
G_{\alpha} (x,y) & \ge & x^\alpha y^\alpha \int ^1{y}\frac{(t-y)^{2\alpha-2} } {t^{2\alpha}} dt \\
&\ge&  x^\alpha y^\alpha \int_{1} ^{\frac{1}{y}}  (1-uy)^{2\alpha-2} du\\
&\ge&  x^\alpha y^\alpha \frac{(1-y)^{2\alpha-1} }{y(2 \alpha -1)} \ge x^{2\alpha-1} \frac{(1-y)^{2\alpha-1} }{2 \alpha -1}.
\end{eqnarray*} 
 C'est \`a dire que nous avons dans ce cas 
 \begin{equation} \label{E1}
  G_{\alpha} (x,y) \ge\frac{ \left(\min (x,y)\right)^{2\alpha-1}  \left (1 - \max (x,y)\right)^{2\alpha-1}}{2\alpha-1}. 
 \end{equation}
 En utilisant la remarque $t \le 1-t $ si et seulement si $ t\le \demi$ nous pouvons crire,  
\begin{equation} \label{E2}
\mathrm{si} \;  \demi \ge \max (x,y) >0 \; \mathrm{alors}\;  
 G_\alpha (x,y) \ge \frac{x ^{2\alpha -1} y ^{2\alpha-1}}{2\alpha-1}
 \end{equation}

\begin{equation} \label{E3}
\mathrm{si}\; 1>\min (x,y) \ge \demi \; \mathrm {alors}\;
 G_\alpha (x,y) \ge \frac{(1-x) ^{2\alpha -1} (1-y) ^{2\alpha-1}}{2\alpha-1}. 
 \end{equation}

La dmonstration du thorme \ref{valeurpropre} nous dit d'autre part 
que 
\begin{align}
&\left( \trace \left( T_N(\varphi_\alpha) \right)^{-1} \right)^s  = \nonumber \\
& \frac{N^{2\alpha s -1}} {( \Gamma ^2 (\alpha ) f_1(1) )^s} 
 \int_0^1 G_\alpha (t_1,t_2) \int _0^1 G_\alpha (t_2,t_3)\cdots
 \int_0^1  G_\alpha (t_{s-1},t_s) G_\alpha (t_s,t_1) dt_s \cdots dt_1.
 \nonumber
 \end{align}
 En d\'ecomposant chacune des int\'egrales ci dessus en une somme de deux int\'egrales,
  l'une o\`u la variable appartient \`a l'intervalle $[0, \demi]$ l'autre o\`u la variable appartient 
  \`a l'intervalle $[\demi, 1]$ on peut \'ecrire  $\left( \trace \left( T_N(\varphi_\alpha) \right)^{-1} \right)^s$
   comme la somme de $2^s$ termes de la forme 
   $$\int_{I_{1}} G_\alpha (t_1,t_2) \int _{I_{2}} G_\alpha (t_2,t_3)\cdots
 \int_{I_{s}}  G_\alpha (t_{s-1},t_s) G_\alpha (t_s,t_1) dt_s \cdots dt_1$$
 o\`u $I_{j}$ est ou \'egal \`a l'intervalle $[0 \demi]$ ou \'egal \`a l'intervalle $ [\demi, 1]$. 
  Alors les minorations \ref{E1}, \ref{E2}, \ref{E3}, permettent d'obtenir   
   $$ \left(\trace\left(\left( T_N(\varphi_\alpha) \right)^{-1} \right)\right)^s
  \ge \frac{2^{s}}{(2\alpha-1)^s} \frac{1}{ \left( \Gamma ^2 (\alpha) f_{1}(1)\right)^s}
 \frac{1}{\left(2^{4\alpha-1} (4\alpha-1)\right)^s} N^{2\alpha s},$$
 ce qui nous donne la minoration 
 $$  \Lambda_{\alpha,N} \ge \frac{2}{2\alpha-1} \frac{1}{ \Gamma ^2 (\alpha) f_{1}(1)}
 \frac{1}{2^{4\alpha-1} (4\alpha-1)} N^{2\alpha}.$$

 \section{D\'emonstration du thoreme \ref{asymptotique}}
 En fait le thorme \ref{asymptotique} a \'et\'e d\'emontr\'ee pour les $\alpha \in \mathbb N$ tendant vers l'infini 
 par Widom et B\"{o}ttcher dans l'article \cite{BoW}. Lorsque cet article a \'et\'e \'ecrit le noyau 
 $G_{\alpha}$ n'\'etait \'etabli que pour $\alpha \in \mathbb N$. N\'eanmoins  nous allons pouvoir r\'eutiliser ici 
 certains des arguments  de cette preuve qui sont toujours valables si on imagine que $\alpha$ est un r\'eel qui tend vers l'infini et non plus
 un entier.\\
 Nous avons d\'ej\`a vu que
 \begin{equation} \label{libel} 
\tilde G_{\alpha} (x,y) = \frac{1}{2} 
G_{\alpha}\left( \frac{1+x}{2}, \frac{1+y}{2}\right) = \frac{1}{4^{2\alpha}} \frac{1}{\Gamma^2 (\alpha)} H_{\alpha}(x,y).
  \end{equation}
  avec 
  \begin{equation} \label{libellule}
  H_{\alpha}(x,y) = (1+x)^\alpha (1+y)^\alpha \int_{\max(x,y)}^1 \frac{ (t-x)^{\alpha-1} (t-y)^{\alpha-1} }{(1+t)/2)^{2\alpha}}dt.
  \end{equation}
 
 En calculant la drive logarithmique de la fonction  
 $t\rightarrow \frac{(t-x) (t-y)}{\left((1+t)/2\right)^2} $
 on montre que la fonction 
 $$(x,y,t) \rightarrow  \frac{(t-x) (t-y)}{\left((1+t)/2\right)^2} $$
 atteint son maximum au seul point $t=1$, $x=0$, $y=0$. C'est \`a dire que si $\delta$ est un rel fix entre $0$ et $1$ on sait qu'en dehors d'un ensemble $\vert t-1 \vert \le \epsilon$, 
 $ \vert x \vert \le \epsilon$, $ \vert y \vert \le \epsilon$, o\`u $0<\epsilon< \demi$ un r\'eel bien choisi on a 
 $$  \frac{(t-x) (t-y)}{\left((1+t)/2\right)^2} < 1- \delta. $$
Puisque $\alpha$ tend vers plus l'infini on peut supposer $\alpha>1$ et donc 
\'ecrire 
$$ \left (\frac{(t-x) (t-y)}{\left((1+t)/2\right)^2} \right)^{\alpha-1}< {1- \delta}^{\alpha-1} $$
ce qui donne 
$$
H_{\alpha}(x,y) =(1+x)^\alpha (1+y)^\alpha  1_{\epsilon}(x) 1_{\epsilon}(y) 
\int_{1-\epsilon}^1 \frac{(t-x)^{\alpha-1} (t-y)^{\alpha-1} }{\left((1+t)/2\right) ^{2\alpha}} dt +O\left( (1-\delta )^\alpha\right),$$
o\`u $1_{\epsilon}$ est la fonction caract\'eristique de l'intervalle $[-\epsilon, \epsilon]$. Avec le changement de variable
$ t = 1-\tau$ on obtient 
\begin{align} \label{jojo}
H_{\alpha}(x,y) & = (1-x^2)^\alpha (1-y^2)^\alpha  1_{\epsilon}(x) 1_{\epsilon}(y) \times \\
& \times \int_{0}^\epsilon \Bigl( \left( 1- \frac{\tau}{1-x} \right) \left( 1- \frac{\tau}{1-y} \right) /  
 \left( 1-\frac{\tau}{2} \right)^2 \Bigr)^\alpha \frac{ d\tau}{ \left( 1- \frac{\tau}{1-x} \right) \left( 1- \frac{\tau}{1-y} \right)}
 +O\left( (1-\delta )^\alpha\right). \nonumber
 \end{align}
 L'objectif est alors de mettre la partie principale de$H_{\alpha}$ sous la forme $f g $ et d'utiliser une 
 propri\'et\'e identique \`a la formule \ref{produit}
 En reprenant les calculs de \cite{BoW} on obtient 
\begin{equation}\label{decomposition}
  H_{\alpha}(x,y) = \frac{1}{\alpha} (1-x^2) ^\alpha (1-y^2)^\alpha \left(1+O(x)+O(y)\right) 
 +O(\frac{1}{\alpha^2})
 \end{equation}
 toujours uniform\'ement en $x$ et $y$ pour $ (x,y) \in [-1,1] \times [-1,1]$. \\
Si  $M_{2}$ et $M_{3}$, d\'esignent les op\'erateurs int\'egraux de noyaux respectif
 $$  O(x) (1-x^2)^\alpha (1-y^2)^\alpha, \quad O(y) (1-x^2)^\alpha (1-y^2)^\alpha.$$
  On a 
 $$  \Vert M_{2}\Vert^2 = \int _{-1}^1O(x^2) (1-x^2)^{2\alpha} dx \int _{-1}^1 (1-y^2)^{2\alpha} dy =
 O\left( \frac{1}{4\alpha} \sqrt {\frac{\pi} {2\alpha} }  \right) \sqrt {\frac{\pi} {2\alpha} } = O\left( \frac{1}{\alpha^2}\right).$$
 On obtient de m\^{e}me $\Vert M_{3}\Vert^2 = O\left( \frac{1}{\alpha^2}\right).$
 Nous pouvons alors \'ecrire, comme dans la d\'emonstration du corollaire \ref{encadrements}
 \begin{align} \label{majordome}
\trace\left(\left( T_N(\varphi_\alpha) \right)^{-1} \right)^s &=  
\frac{N^{2\alpha s}}{ \left( \Gamma^2 (\alpha)f_1(1)\right)^s} (\frac{1}{4^{2\alpha}} )^s
   \int _{-1}^1 \int _{-1}^1 H_\alpha ( \frac{1+x_1}{2},\frac{1+x_2}{2})\times \nonumber  \\
&\times \int_{-1}^1 H_\alpha (\frac{1+x_2}{2},\frac{1+x_3}{2})
\cdots \int _{-1}^1H_\alpha (\frac{1+x_{s-1}}{2},\frac{1+x_s}{2}) \times  \\
&\times \int_{-1}^1 H_\alpha (\frac{1+x_s}{2},\frac{1+x_1}{2})
dx_1 \cdots dx_s. \nonumber
\end{align}
En utilisant \ref{decomposition} et les r\'esultats 
$$ \Vert M_{2}\Vert^2 =\Vert M_{3}\Vert^2 = O\left( \frac{1}{\alpha^2}\right).$$
On peut \'ecrire, \`a partir de \ref{majordome}
 \begin{align} \label{majordome2}
\trace\left(\left( T_N(\varphi_\alpha) \right)^{-1} \right)^s &=  
\frac{N^{2\alpha s}}{ \left( \Gamma^2 (\alpha)f_1(1)\right)^s} (\frac{1}{4^{2\alpha}} )^s \frac{1}{\alpha^s}
 \left(  \int _{-1}^1 \int _{-1}^1 (1-x_1)^\alpha\ (1-x_2) ^\alpha\times \right.\nonumber  \\
&\times \int_{-1}^1  (1-x_2)^\alpha (1-x_3)^\alpha
\cdots \int _{-1}^1 (1-x_{s-1})^\alpha,(1-x_s)^\alpha \times  \\
&\times \left.\int_{-1}^1 (1-x_s)^\alpha,(1-x_1)^\alpha 
dx_1 \cdots dx_s. \right) \times \nonumber\\
& \times \left( 1+O(1) \right) \nonumber
\end{align}
 Nous avons d'autre part 
$$ \int_{-1}^1 (1-x^2) ^{2\alpha} dx = \sqrt { \frac{\pi}{2\alpha}} \left( 1 + O(\frac{1}{\alpha}) \right).$$
 L'\'egalit\'e \ref{majordome2} permet donc d'\'ecrire, avec toujours la m\^eme notation pour la plus grande valeur propre 
  $$ \Lambda_{\alpha,N} = \frac {N^{2\alpha}}{f_1(1)}
  \frac{1}{ \alpha \Gamma^2(\alpha) 4^{2\alpha}} \left( \sqrt{\frac{\pi}{2\alpha}}+ O\left(\frac{1}{\alpha}\right)\right),$$
 ou encore 
  $$  \Lambda_{\alpha,N}= \frac {N^{2\alpha}}{f_1(1)}\left( \frac{e}{4\alpha} \right) ^{2\alpha} \frac{1}{\sqrt{8\pi \alpha}} \left( 1+O\left( \frac{1}{\sqrt \alpha} \right) \right). 
$$

   \bibliography {Toeplitzdeux}

\begin{thebibliography}{10}

\bibitem{Bot}
A.~B\"ottcher.
\newblock The constants in the asymptotic formulas by {R}ambour and {S}eghier
  for the inverse of {T}oeplitz matrices.
\newblock {\em Integr. equ. oper. theory}, 99:43--45, 2004.

\bibitem{BoG}
A.~B\"ottcher and J~Unterberger S.~Grudsky, E. A.~Maksimenko.
\newblock The first order asymptotics of the extreme eigenvectors of certain
  hermitian {T}oeplitz matrices.
\newblock {\em Integr. equ. oper. theory}, 63:165--180, 2009.

\bibitem{Bo.3}
A.~B\"ottcher and B.~Silbermann.
\newblock {\em Introduction to large truncated matrices}.
\newblock Springer Verlag, 1999.

\bibitem{BoW}
A.~B\"ottcher and H.~Widom.
\newblock From {T}oeplitz eigenvalues through {G}reen's kernels to higher-order
  {W}irtinger-{S}obolev inequalities.
\newblock {\em Oper. Th. Adv. Appl.}, 171:73--87, 2006.

\bibitem{BoW2}
A.~B\"ottcher and H.~Widom.
\newblock On the eigenvalues of certain canonical higher-order ordinary
  differential operators.
\newblock {\em J. Math. Anal. Appl.}, 322:990--1000, 2006.

\bibitem{CFL}
R.~Courant, K.~Friedrichs, and H.~Lewy.
\newblock {\"U}ber die partiellen {D}ifferenzengleichungen der mathematischen
  {P}hysik.
\newblock {\em Math. Ann.}, 100:32--74, 1928.

\bibitem{GoSe}
I.~Gohberg and A.~A. Semencul.
\newblock The inversion of finite {T}oeplitz matrices and their continual
  analogues.
\newblock {\em Matem. Issled.}, 7:201--233, 1972.

\bibitem{GS}
U.~Grenander and G.~Szeg{\"o}.
\newblock {\em {T}oeplitz forms and their applications}.
\newblock Chelsea, New York, 2nd ed. edition, 1984.

\bibitem{KMS}
M.~Kac, W.~L. Murdoch, and G.~Szeg{\"o}.
\newblock On the eigenvalues of certain hermitian forms.
\newblock {\em J. rat. Mech. Analysis}, 2:767--800, 1953.

\bibitem{KaRS}
D.~Kateb, P.~Rambour, and A.~Seghier.
\newblock Asymptotic behavior of the predictor polynomial associated to regular
  symbols.
\newblock Pr\'epublications de l'Universit\'e Paris-Sud, 2003.

\bibitem{Ld}
H.J. Landau.
\newblock Maximum entropy and the moment problem.
\newblock {\em Bulletin ({N}ew {S}eries) of the american mathematical society},
  16(1):47--77, 1987.

\bibitem{ML3}
A.~Martinez-Finkelshtein, K.~T.~R McLaughlin, and E.~B. Saff.
\newblock Asymptotics of orthogonal polynomials with respect to an analytic
  weight with algebraic singularities on the circle.
\newblock {\em Internat. Math. Research Notices}, 2006.

\bibitem{Part}
S.~Parter.
\newblock Extreme eigenvalues of {T}oeplitz forms and applications to elliptic
  difference equations.
\newblock {\em Trans. Amer. Math. Soc.}, 99:153--192, 1961.

\bibitem{Part3}
S.~Parter.
\newblock On the extreme eigenvalues of {T}oeplitz matrices.
\newblock {\em Trans. Amer. Math. Soc.}, 100:263--270, 1961.

\bibitem{Part2}
S.~Parter.
\newblock On the extreme eigenvalues of truncated {T}oeplitz matrices.
\newblock {\em Bull. Amer. Math. Soc.}, 67:191--196, 1961.

\bibitem{RS04}
P.~Rambour and A.~Seghier.
\newblock Formulas for the inverses of {T}oeplitz matrices with polynomially
  singular symbols.
\newblock {\em Integr. equ. oper. theory}, 50:83--114, 2004.

\bibitem{RS10}
P.~Rambour and A.~{S}eghier.
\newblock Inverse asymptotique des matrices de {T}oeplitz de symbole $(1-\cos
  \theta)^\alpha f_{1},$ $ \frac{-1}{2}<\alpha\le \frac{1}{2}$, et noyaux
  int\'egraux.
\newblock {\em Bull. des {S}ci math,}, 134:155--188, 2008.

\bibitem{RS12}
P.~Rambour and A.~Seghier.
\newblock Asymptotic inversion of toeplitz matrices with one singularity in the
  symbol.
\newblock {\em C. R. Acad. Sci. Paris}, 347, ser. I:489--494, 2009.

\bibitem{SpSt}
F.~L. Spitzer and C.~J. Stone.
\newblock A class of {T}oeplitz forms and their applications to probability
  theory.
\newblock {\em Illinois J. Math.}, 4:253--277, 1960.

\bibitem{GZ}
G.~Szeg\"{o}.
\newblock A {T}oeplitz fle form\'{a}kr\'{o}l.
\newblock {\em Mathematikai s termszettudom\'anyi ertesit\"o}, 35:185--222,
  1917.

\bibitem{Wid}
H.~Widom.
\newblock On the eigenvalues of certain hermitian operators.
\newblock {\em Trans. Amer. Math. Soc.}, 88:491--522, 1958.

\bibitem{W3}
H.~Widom.
\newblock Extreme eigenvalues of translation kernels.
\newblock {\em Trans. amer. Math. Soc.}, 100:252--262, 1961.

\bibitem{W2}
H.~Widom.
\newblock Extreme eigenvalues of {N}-dimensional convolution operators.
\newblock {\em Trans. Amer. Math. Soc.}, 106:391--414, 1963.

\end{thebibliography}
  \end{document}